\documentclass[a4paper,fleqn, final]{cas-dc}

\usepackage[authoryear,longnamesfirst]{natbib}
\usepackage{tikz}
\usepackage{xcolor}
\usepackage{array, colortbl, xcolor, booktabs}
\usepackage{pgfplots}
\usepackage{afterpage}
\usepackage{amsmath}
\usepackage{float}
\usepackage{caption}
\usepackage{booktabs}
\usepackage{multirow}
\usepackage{algorithm}
\usepackage{algorithmic}
\usepackage{graphicx}
\usepackage{caption}
\usepackage{subcaption}
\usepackage{threeparttable}

\usepackage{color}
\definecolor{Yingxiang}{rgb}{1, 0, 0}
\definecolor{other}{rgb}{0, 0, 1}
\definecolor{zhangshuai}{rgb}{1, 0, 1}

\begin{document}

\let\WriteBookmarks\relax
\def\floatpagepagefraction{1}
\def\textpagefraction{.001}

\shorttitle{Neural Networks for Solving SIEPs}    

\shortauthors{S. Zhang and Y. Xu}  

\title [mode = title]{Orthogonal Constrained Neural Networks for Solving Structured Inverse Eigenvalue Problems}

\author[1]{Shuai Zhang}

\ead{zhangs216@nenu.edu.cn}

\author[1]{Xuelian Jiang}
\ead{jiangxl133@nenu.edu.cn}

\author[1]{Hao Qian}
\ead{qianhao123@nenu.edu.cn}

\author[1]{Yingxiang Xu}
\ead{yxxu@nenu.edu.cn}
\cormark[1]
\cortext[1]{Corresponding author. School of Mathematics and Statistics, Northeast Normal University, Changchun 130024, P. R. China}

\affiliation[1]{organization={School of Mathematics and Statistics, Northeast Normal University},
                city={Changchun},
                postcode={130024}, 
                country={China}}



\begin{abstract}
      This paper introduces a novel neural network for efficiently solving Structured Inverse Eigenvalue Problems (SIEPs). The main contributions lie in two aspects: firstly, a unified framework is proposed that can handle various SIEPs instances. Particularly, an innovative method for handling nonnegativity constraints is devised using the ReLU function. Secondly, a novel neural network based on multilayer perceptrons, utilizing the Stiefel layer, is designed to efficiently solve SIEP. By incorporating the Stiefel layer through matrix orthogonal decomposition, the orthogonality of similarity transformations is ensured, leading to accurate solutions for SIEPs. Hence, we name this new network Stiefel Multilayer Perceptron (SMLP). Furthermore, SMLP is an unsupervised learning approach with a lightweight structure that is easy to train. Several numerical tests from literature and engineering domains demonstrate the efficiency of SMLP.
\end{abstract}
\maketitle
\section{Introduction}
A Structured Inverse Eigenvalue Problem (SIEP) consists of reconstructing a structured matrix from prescribed spectral data.
The SIEP is a significant research topic with a wide range of applications such as inverse Sturm-Liouville problems \citep{hald1972discrete,hald1989solution}, graph theory \citep{hogben2005spectral}, engineering structural design \citep{joseph1992inverse,collins2016design}, geophysics \citep{parker1981numerical}, compressed sensing \citep{baraniuk2017compressive}, and microwave filter network transformation \citep{lamecki2004fast,kozakowski2005eigenvalue,pfluger2015coupling}.
The objective of the SIEP is to construct a matrix that maintains both the specific structure as well as the given spectral property.
The theoretical study focuses on the necessary and sufficient conditions for the existence of solutions, under the impact of the entry values, and locations of prescribed entries in a structured matrix. 
The numerical study, which aims to complete such a construction assuming the existence of solutions, is complex and demanding because of the complex pattern in the structure.

\begin{figure}[h]
       \centering
       \includegraphics[width=0.5\textwidth]{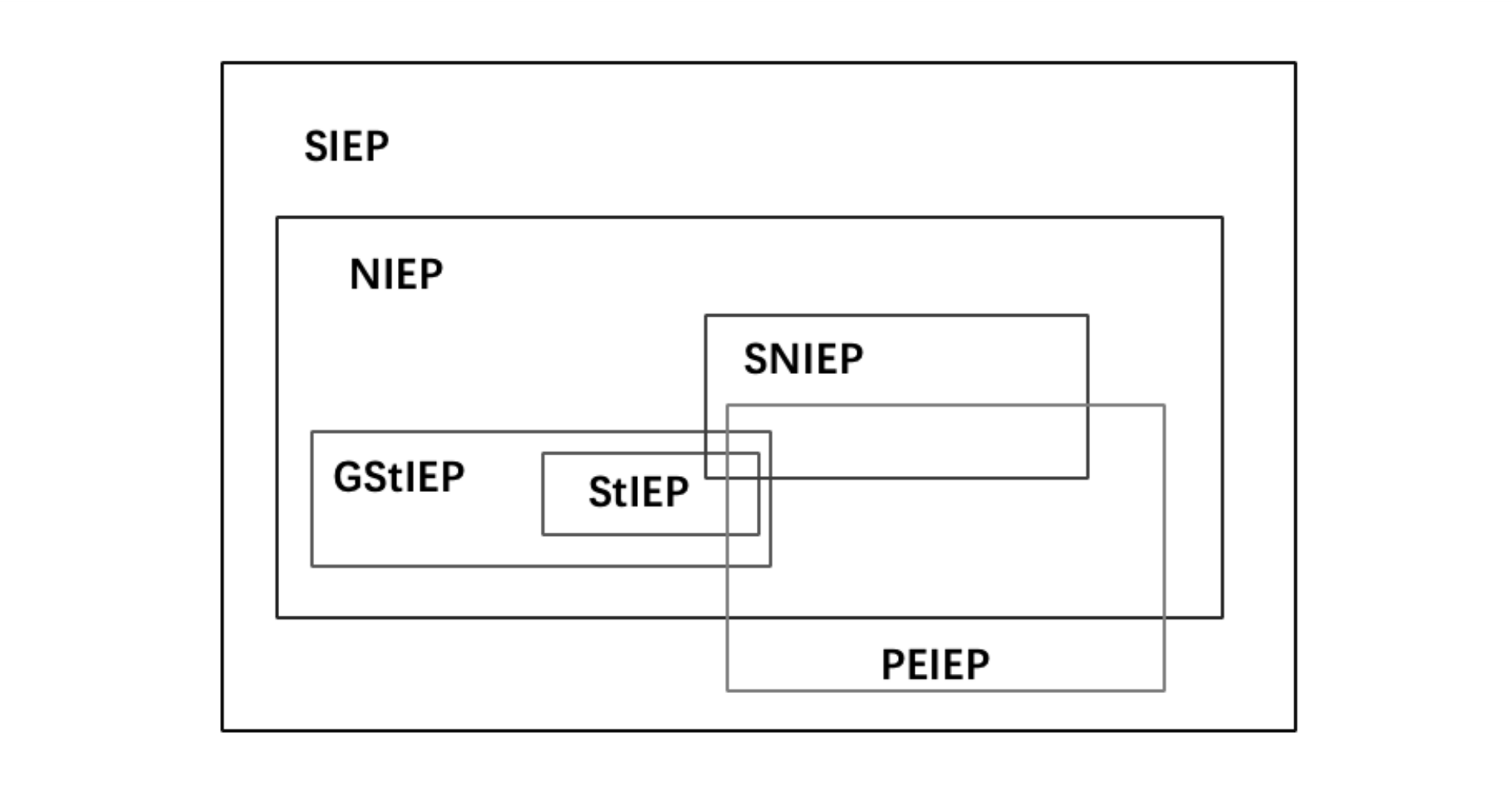} 
       \caption{Classification of SIEPs under consideration.}\label{CLASSIEP}
\end{figure}

For a detailed classification of the SIEP, we refer the reader to  \citet{chu2002structured,chu2001inverse}.
In this paper, we focus on SIEP with prescribed entries (PEIEP), nonnegative inverse eigenvalue problems (NIEP), and nonnegative inverse eigenvalue problems with prescribed entries (NPEIEP). Specifically, we explore symmetric nonnegative inverse eigenvalue problems (SNIEP), stochastic inverse eigenvalue problems (StIEP) and generalized stochastic inverse eigenvalue problems (GStIEP).
For stochastic matrices, the inverse eigenvalue problem is particularly difficult, as evidenced by the complexities involved in the well-known result on existence by  \citet{karpelevich1951characteristic}.
The diversity of structural constraints account for the various types of SIEPs elucidated in this paper, as illustrated in Figure \ref{CLASSIEP}.
However, the loss functions corresponding to different type of SIEPs vary when solving from an optimization perspective, meaning that a specific loss function is tailored to particular type of SIEPs.

Neural networks have gained significant attention due to their universal approximation capabilities and powerful expression abilities. 
Neural networks are widely used in fields such as computer vision, natural language processing, and image processing. 
Recently, advancements have been made in developing neural networks for solving eigenvalue problems. 
For example, a functional neural network computing some eigenvalues and eigenvectors of a special real matrix \citep{liu2005functional}, neural networks computing a generalized eigenvalue problem \citep{hara2012separation}, neural networks based on the power and inverse power methods for solving linear eigenvalue problems \citep{yang2023neural}, and adversarial neural networks are employed to solve topology optimization of eigenvalue problems \citep{hu2024adversarial}. 
To prevent the gradient explosion/vanishing issues often occurred in network training, the orthogonal constraints on network parameters are proven to be good strategies \citep{mhammedi2017efficient,huang2020controllable,zhang2021orthogonality}.
The neural networks for solving inverse eigenvalue problems associated to partial differential equations are detailed in Section 2.
To our best knowledge, up to now there is no research on neural networks solver for algebraic inverse eigenvalue problems, especially for the SIEP. 

This paper presents, for the first time, a well-designed neural network for solving the SIEP. We first transform the SIEP into an optimization problem on the Stiefel manifold and propose a more generic loss function in a unified form, which allows us to address all SIEPs  under consideration. Different from a soft constraint approach where the constraints are incorporated in the loss function, we use the idea of the hard constraint by embedding the orthogonal constraints directly into the structure of the network through matrix orthogonal decomposition techniques.
This innovation leads the output of the network to strictly satisfying the orthogonal constraints, that is to say, it belongs to the Stiefel manifold. 
Additionally, this is an unsupervised learning process that does not rely on labeled datasets for training. 
As it is well known, the SIEP is hard to optimize due to their nonconvex property. While our approach exhibits superior capabilities in efficiently finding the global optimizer for the SIEP, which shows great potential to extending to other type of matrix inverse eigenvalue problems.

This paper is outlined of as follows. We present related works in Section 2. Section 3 provides the necessary background. Section 4 details the design of the unified loss function and proposes the new neural network SMLP with Stiefel Layer. 
Section 5, to demonstrate the efficiency of the proposed SMLP, presents and analyzes the experimental results on solving several typical SIEPs from literature and an engineering application from microwave filter design.
Finally, we draw conclusions in Section 6.

\section{Related works}

\paragraph{The solvability of the IEP.} The fundamental issues related to the matrix inverse eigenvalue problem mainly include two aspects: the theory of solvability and the practice of computability. The analysis of solvability theory aims to identify the necessary or sufficient conditions for the existence of solutions to the IEP. 
Nonnegative matrices have attracted significant attention due to their important applications in various fields such as game theory, Markov chains, and economics. 
In 1949, Suleimanova extended Kolmogorov's problem and proposed the second major conjecture in matrix theory: to find a set of necessary and sufficient conditions such that a given closed conjugate set of complex numbers is the spectrum of a nonnegative matrix \citep{suleimanova1949stochastic}. 
Subsequently, the NIEP has garnered considerable scholarly attention.
Sufficient conditions for the existence of solutions to the NIEP were explored \citep{kellogg1971matrices,borobia1995nonnegative,soto2003existence,marijuan2014sufficient,marijuan2007map}, while \citet{loewy1978note} studied some necessary conditions. 
\citet{swift1972location} and \citet{shang2011row} proposed sufficient conditions for the solvability of the IEP for stochastic matrices. 
Some necessary and sufficient conditions for the IEP of distance matrices have also been investigated \citep{jaklivc2012note,hayden1999methods,nazari2014inverse}.

\paragraph{The computability of the IEP.} The objective of the feasibility analysis is to develop an algorithm that constructs a target matrix meeting spectral constraints for provided solvable spectral data using a numerically stable approach.
One category involves constructive algorithms \citep{soto2006applications,laffey2007construction,soto2003existence}, which, under certain sufficient conditions, can construct nonnegative matrices with desired spectra. 
However, these constructive methods are only applicable for matrices with specific structures.
Another category involves optimization algorithms, which transform the IEP into an optimization problem and then solve it using optimization techniques within the Riemannian manifold framework. It  includes Riemannian Gradient Descent \citep{huang2018blind} and Riemannian Conjugate Gradient methods \citep{sato2016dai,yao2019riemannian}.
Both above methods are geometric optimization methods. As a drawback of these geometric methods, it requires to compute analytical solutions for geometric quantities, such as the gradient direction and geodesics defined by the Riemannian metric \citep{kanamori2014numerical}.

\paragraph{Establishment of the loss for the IEP.}
\citet{chu1990projected} proposed the IEP of symmetric matrices under spectral constraints. 
The given spectral data is used as the entries of a diagonal matrix to construct an isospectral manifold. 
The loss function is defined as the distance between the isospectral manifold and the set of all symmetric matrices. 
\citet{zhao2018riemannian} studied the IEP of finding a nonnegative matrix from a given realizable spectrum. 
They defined the IEP as solving an underdetermined constrained nonlinear matrix equation on a matrix manifold, and constructed a loss function based on this equation.
A special type IEP that the target matrix has specified entries while satisfying spectral constraints was investigated by  
\citet{chu2004gradient}, where the problem is transformed into minimizing the distance between an isospectral matrix with preset eigenvalues and an affine matrix with given elements. 
\citet{chen2011isospectral} further extended the work of \citet{chu2004gradient}, studying the IEP of symmetric nonnegative matrices with preset elements, ensuring nonnegativity through the hadamard product of matrices.

\paragraph{Neural networks for the IEP in PDEs.}
\citet{ossandon2016neural} proposed a numerical method utilizing an artificial neural network to address the inverse problem linked to computing the Dirichlet eigenvalues of the anisotropic Laplace operator. 
By repeatedly solving the direct problem to obtain a set of predefined eigenvalues, a neural network was designed to find the appropriate components of the anisotropic matrix related to the Laplace operator, thereby solving the associated inverse problem.
\citet{pallikarakis2024application} investigated the numerical solution of inverse eigenvalue problems from the perspective of machine learning. 
It focuses on two distinct issues: the inverse Sturm-Liouville eigenvalue problem with symmetric potentials and the inverse transmission eigenvalue problem with spherically symmetric refractive indices. 
The study employs supervised regression models such as k-Nearest Neighbors, Random Forests, and MultiLayer Perceptron.
However, no attempt on solving algebraic inverse eigenvalue problems using neural networks published.

\section{Preliminaries}
To start with, we will first briefly introduce some fundamental concepts used in this paper, including the basic definitions in matrix theory, the matrix decomposition techniques, the definition of the Stiefel manifold. 
Additionally, we will explore the fundamental properties of neural networks to lay the groundwork for subsequent investigations.

\subsection{Classification of matrices}
Several type of matrices considered in this paper are defined below, see for example \citep{horn2012matrix}. 

\paragraph{Nonnegative  matrix.} A nonnegative matrix is defined as a matrix where all entries are nonnegative.

\paragraph{Stochastic matrix.} A matrix $A$ is stochastic if all of its entries are nonnegative, and the entries of each row sum to 1.
       Let $A$ be a stochastic matrix, then 1 is an eigenvalue of $A$. If $\lambda$ is an eigenvalue of $A$, then $|\lambda| \leq 1$.

\paragraph{Generalized stochastic matrix.} A generalized stochastic matrix is a matrix whose elements are nonnegative, and the sum of the elements in each row is a constant number.

\paragraph{Euclidean distance matrix (EDM).} A matrix \(D = (d_{ij}) \in \mathbb{R}^{n \times n}(1\leq i,j\leq n)\) is a Euclidean distance matrix if it represents the squared Euclidean distances \(d_{ij} = \|x_i - x_j\|^2\) between \(n\) points \(x_1, x_2, \ldots, x_n \in \mathbb{R}^r\), where \(\|\cdot\|\) denotes the Euclidean norm.
       \(D\) is a nonnegative and symmetric matrix, and has a zero main diagonal, implying the sum of its eigenvalues is zero.

\subsection{Matrix decomposition}
Among the various methods for matrix decomposition, we review in this paper the QR decomposition, Singular Value decomposition and Real Schur decomposition.

\paragraph{QR decomposition.}
Let $A \in \mathbb{R}^{m \times n}$ with $m \geq n$ be a rectangular matrix. Then \( A \) has a QR decomposition:
\begin{equation}\label{QR}
A = Q \begin{bmatrix}
R \\
0
\end{bmatrix},
\end{equation}
where \( Q \in \mathbb{R}^{m \times m} \) is an orthogonal matrix and \( R \in \mathbb{R}^{n \times n} \) is an upper triangular matrix with nonnegative diagonal elements. Furthermore, when \( m = n \) and \( A \) is nonsingular, the decomposition is unique.

\paragraph{Singular value decomposition (SVD).}
Let $A \in \mathbb{R}^{m \times n}$ with $m \geq n$. Then there exist orthogonal matrices $U \in \mathbb{R}^{m \times m}$ and $V \in \mathbb{R}^{n \times n}$ and a diagonal matrix $\Sigma=\operatorname{diag}\left(\sigma_1, \cdots, \sigma_n\right) \in \mathbb{R}^{m \times n}$ with $\sigma_1 \geq \sigma_2 \geq \cdots \geq \sigma_n \geq 0$, such that
\begin{equation}\label{SVD}
A=U \Sigma V^T
\end{equation}
holds. The column vectors of $U=\left[{u}_1, \cdots, {u}_m\right]$ are called the left singular vectors and similarly $V=\left[{v}_1, \cdots, {v}_n\right]$ are the right singular vectors. The values $\sigma_i$ are called the singular values of $A$. 

\paragraph{Real Schur decomposition.}
Suppose \( A \in \mathbb{R}^{n \times n} \). Then there exists an orthogonal matrix \( Q \in \mathbb{R}^{n \times n} \) such that
\begin{equation}\label{Sde}
Q^{\mathrm{T}} A Q = \left[\begin{array}{cccc}
R_{11} & R_{12} & \cdots & R_{1n} \\
& R_{22} & \cdots & R_{2n} \\
& & \ddots & \vdots \\
& & & R_{nn}
\end{array}\right],
\end{equation}
where \( R_{ij} \) ($1 \leq i,j\leq n$) are real numbers. In particular, if \( A \) is symmetric (i.e., \( A = A^{\mathrm{T}} \)), then \( Q^{\mathrm{T}} A Q \) is orthogonally similar to a diagonal matrix.
Real Schur decomposition indicates that any real matrix can be orthogonally transformed into an upper triangular matrix.

\subsection{Neural network}
To introduce the foundational aspects of neural networks, we consider the Multilayer Perceptron (MLP). 
The MLP, features fully connected input, hidden, and output layers, making it one of the most widely used neural network models. 

\paragraph{Multilayer perceptron neural network.}
\begin{figure}[h]
       \centering
       \includegraphics[width=0.4\textwidth]{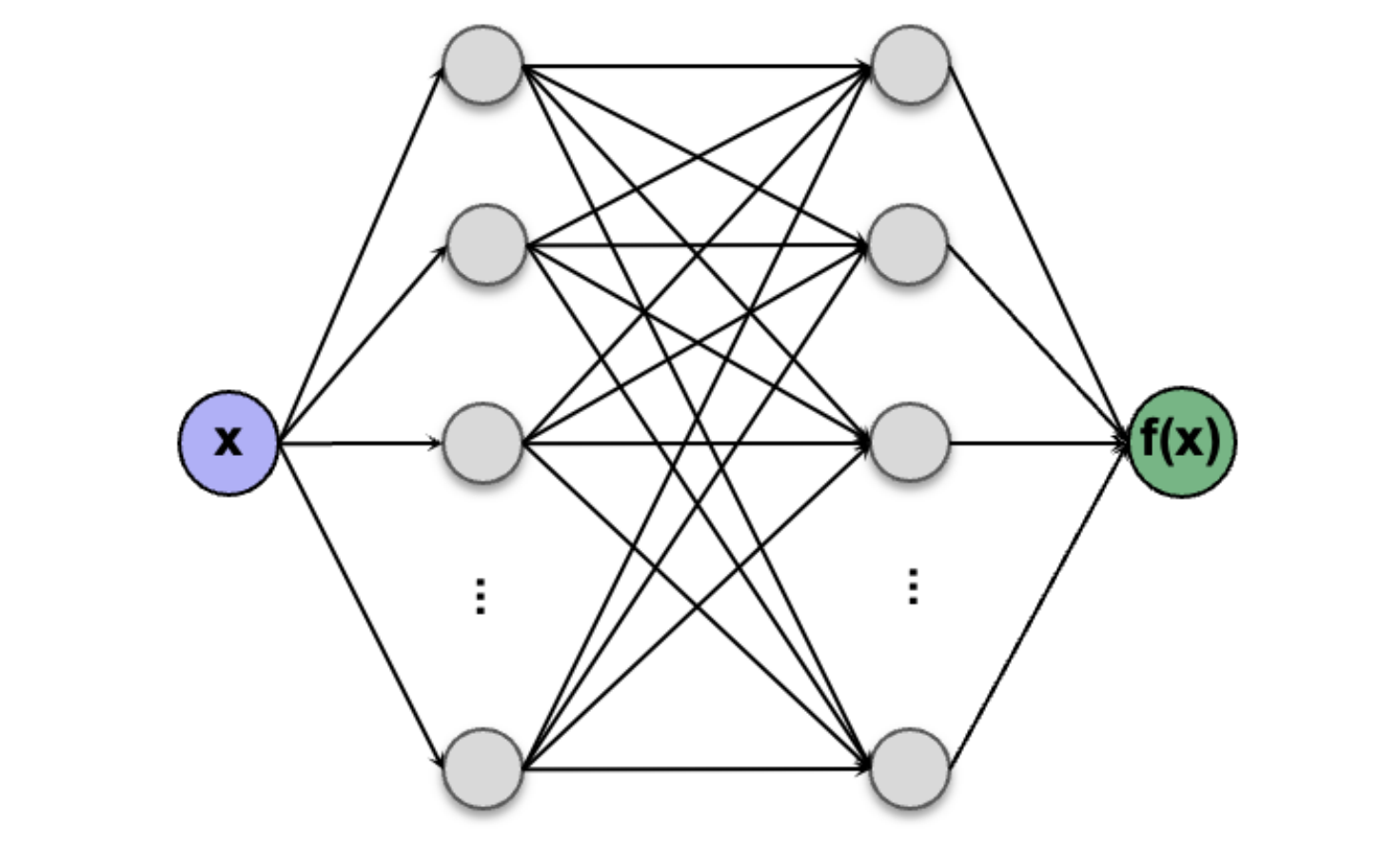} 
       \caption{The structure of the MLP consists of a fully connected network composed of an input layer, hidden layers, and an output layer.}\label{MLP}
\end{figure}
Its structure is divided into three types of layers that are fully connected among themselves: the input layer, where the model inputs are introduced; the output layer, where the results of the trained model are obtained; and the hidden layers, which are the intermediate layers between the preceding ones, as shown in Figure \ref{MLP}.
The MLP, denoted as \(\mathcal{F} = \{f: \mathbb{R}^p \rightarrow \mathbb{R}^d\}\), is represented by:
$$
f(x) = W_{\mathcal{L}}  *  W_{\mathcal{L}-1} * \Phi * \cdots \Phi * W_1 * \Phi * W_0(x),
$$
where \(x = [x_1, \ldots, x_p]\in \mathbb{R}^p\), \(W_i(x) = \omega_i x + b_i\) is an affine transformation; \(\omega_i \in \mathbb{R}^{k_{i+1} \times k_i}\) is the weight matrix; \(b_i \in \mathbb{R}^{k_{i+1}}\) is the bias vector for \(i = 0, 1, \ldots, \mathcal{L}\); and \(\Phi\) is the activation function.
The input and output layers have widths \(k_0 = p\) and \(k_{\mathcal{L}+1} = d\), respectively. 
\(k_i\) (\(i = 1, 2, \ldots, \mathcal{L}\)) represents the number of neurons in the \(i\)-th hidden layer. \(\mathcal{L}\) denotes the number of hidden layers, also known as the depth of the network. The width of the network refers to the maximum width of the hidden layers, i.e. \(\max \{k_1, k_2, \ldots, k_{\mathcal{L}}\}\).
The total number of neurons in the network is \(\sum_{i=1}^{\mathcal{L}} k_i\), and the total number of parameters equals \(\sum_{i=0}^{\mathcal{L}} k_{i+1} \times (k_i + 1)\).

\paragraph{Activation functions.}
The following introduces two commonly used activation functions. 

A ReLU is mathematically defined as the maximum of zero and the input value :
       $$
       \Phi(x):={\tt ReLU}(x)=\max (0, x) .
       $$
The mathematical expression for ${\tanh}$ is: 
       $$
       \Phi(x) :={\tanh}(x)= \frac{e^x - e^{-x}}{e^x + e^{-x}} .
       $$
For additional activation functions, one can refer to \citet{kunc2024three}.

\subsection{Stiefel manifold}
     The Stiefel manifold, denoted as $\mathcal{S}_{n, r}$, comprises $n \times r$ rectangular matrices whose column vectors are orthonormal, i.e.,
     \begin{equation}\label{Sti}
     \mathcal{S}_{n, r}=\left\{W \in \mathbb{R}^{n \times r}: W^T W=\boldsymbol{I}_r\right\},
     \end{equation}
     where it is assumed that $n \geq r$, \(\boldsymbol{I}_r\) represents the \(r\) dimensional identity matrix and $W$ is a semi-orthogonal matrix when $r<n$. 
     The Stiefel manifold when $n=r$ is equivalent to the orthogonal group, which is the set of all orthogonal matrices, and  $W$ is an orthogonal matrix.
     The optimization problem involving the Stiefel manifold can be expressed as:
     $$
     \min _{W\in \mathcal{S}_{n, r}} g(W),
     $$
     where $g(\cdot)$ is a real-valued loss function defined on the set of rectangular matrices $\mathbb{R}^{n \times r}$.

\section{Methodology}
In this section, we propose a unified framework for solving the SIEP. We discuss in detail the design of specific loss functions under different constraints and propose a unified loss function to handle the SIEP under these constraints, which is one of our main contributions. Another main contribution is that
we devise the MLP with a Stiefel layer (SMLP) to ensure that the network outputs satisfy orthogonality constraints. 
We describe the architecture and training process of SMLP, demonstrating how SLMP can be effectively applied to solve SIEPs.
\subsection{Loss function establishment in SIEPs}
We establish the loss function starting from the SIEP with prescribed entries. 
Further, we apply nonnegativity constraints to the matrix elements and propose restrictions on the sum of elements in each row. 
For these constraints, we design a unified loss function to meet these conditions. 
The nonnegativity constraint  is fulfilled using the ReLU function, a novel approach.

\paragraph{Design of loss functions under various constraints.}
Consider first the inverse eigenvalue problem with prescribed entries (PEIEP). Given the spectrum $\sigma=\{\lambda_1,\ldots,\lambda_n\}$. Let \(\Lambda\) be an upper triangular matrix with \(\lambda_1, \ldots, \lambda_n\) as its main diagonal elements. The isospectral manifold corresponding to $\Lambda$ is defined as:
\begin{equation*}
  \mathcal{M}(\Lambda):=\left\{Q \Lambda Q^T \mid Q \in \mathcal{O}(n)\right\},
\end{equation*}
where $\mathcal{O}(n)$ denotes the set of all orthogonal matrices in $\mathbb{R}^{n \times n}$.
Given the index set $L=\{i_v,j_v\}_{v=1}^l$, where $1 \leq i_v,j_v \leq n$ and $l \leq n^2$, and the corresponding set of values $\omega=\{\omega_1,\ldots,\omega_l\}$.
Define the matrix $\Omega \in \mathbb{R}^{n\times n}$ such that $\Omega_{i_v,j_v}=\omega_v$ for $v=1,\ldots,l$, with all other elements being zero. Let the matrix $S \in \mathbb{R}^{n \times n}$ be such that $S_{i_v,j_v}=0$ for $v=1,\ldots,l$, with all other elements being set to one.
Let the matrix $\Gamma \in \mathbb{R}^{n\times n}$ satisfy $\Gamma+\Omega=E$, where $E$ is the $n \times n$ matrix with all elements equal to one. 

Let
\begin{equation}\label{pi}
  \pi_{L,\omega}=\{\Omega+S\circ M | M\in \mathbb{R}^{n\times n}\},
\end{equation}
where $\circ$ denotes the Hadamard product.
Thus, the PEIEP has a solution if and only if:
\begin{equation*}
  \mathcal{M}(\Lambda) \cap \pi_{L,\omega} \neq \emptyset,
\end{equation*}
that is, the following nonlinear matrix equation has a solution:
\begin{equation}\label{matrix eq}
  Q \Lambda Q^T  -(\Omega+S\circ M ) = 0.
\end{equation}

Then (\ref{matrix eq}) can be reformulated as a problem of minimizing the distance between the isospectral manifold $\mathcal{M}$ and $\pi_{L,\omega}$, expressed as follows:
\begin{equation}\label{minfQM1}
       \begin{aligned}
              \min \quad  & F(Q,M) = \frac{1}{2}\| Q \Lambda Q^T  -(\Omega+S\circ M ) \|_F^2, \\
              \text{s.t.}\quad  &  Q\in \mathcal{O}(n).
       \end{aligned}
\end{equation}
Performing the Schur decomposition, as in Formula (\ref{Sde}), on matrix $M$, we can restate (\ref{minfQM1}) as:
\begin{equation*}
       \begin{aligned}
              \min \quad  &F(Q)  = \frac{1}{2}\| Q \Lambda Q^T  - \Omega - S\circ (Q \Lambda Q^T) \|_F^2, \\
               \text{s.t.}  \quad  &  Q\in \mathcal{O}(n).
       \end{aligned}
\end{equation*}
Since
\begin{equation*}
       \begin{aligned}
       &\| Q \Lambda Q^T  - \Omega - S\circ (Q \Lambda Q^T) \|_F^2 \\
       =&\| Q \Lambda Q^T - S\circ (Q \Lambda Q^T) - \Gamma \circ (Q \Lambda Q^T) - \Omega + \Gamma \circ (Q \Lambda Q^T)\|_F^2 \\
       =&\|(\Omega - \Gamma \circ (Q \Lambda Q^T))\|_F^2,
       \end{aligned}
\end{equation*}
we obtain the following loss function:
\begin{equation}\label{Loss1}
       \begin{aligned}
              \min \quad &F(Q)  = \frac{1}{2}\|\Omega - \Gamma \circ (Q \Lambda Q^T)\|_F^2. \\
       \end{aligned}
\end{equation}
We observe that the ReLU activation function in neural networks has desirable mathematical properties, and thus propose the following loss function for the NPEIEP:
\begin{equation}\label{Loss2}
       \begin{aligned}
         \min \quad & F(Q) = \frac{1}{2}\| S \circ (Q \Lambda Q^T) - \text{ReLU}(S \circ (Q \Lambda Q^T)) \|_F^2 \\
                    &\quad \quad \quad+ \frac{1}{2}\|\Omega - \text{ReLU}(\Gamma \circ (Q \Lambda Q^T))\|_F^2, \\
         \text{s.t.} \quad & Q \in \mathcal{O}(n).
       \end{aligned}
\end{equation}
Next, let $\alpha=[a_1,\ldots,a_n]_{1\times n}$, and $\beta=[b_1,b_2,\ldots,b_n]$, where $b_i$ ($1 \leq i \leq n$) represents the sum of the elements in the $i$-th row of the target matrix. For a  generalized stochastic matrix with a row sum constant $a_i=a$ (which is a stochastic matrix when $a=1$), the following loss function is proposed:
\begin{equation}\label{Loss3}
\begin{aligned} 
  \min \quad &F(Q)  = \frac{1}{2}\| S \circ (Q \Lambda Q^T) - \text{ReLU}(S\circ (Q \Lambda Q^T)) \|_F^2  \\
  &\quad \quad \quad+ \frac{1}{2}\|\Omega - \text{ReLU}(\Gamma \circ (Q \Lambda Q^T))\|_F^2 \\ 
  &\quad \quad \quad + \frac{1}{2}\|\alpha - \beta\|_F^2,\\
  \text { s.t. }   \quad &Q\in \mathcal{O}(n).
\end{aligned}
\end{equation}
\begin{figure*}[h]
       \centering
       \includegraphics[width=1.0\textwidth]{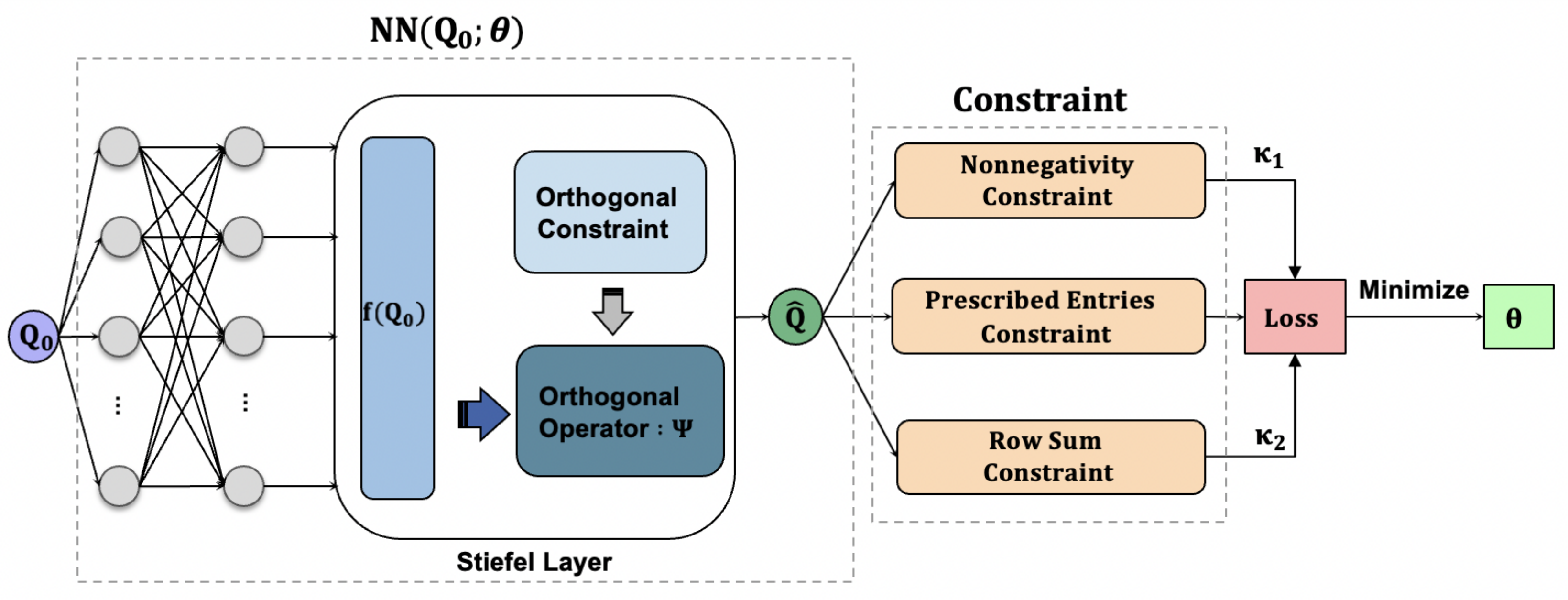} 
       \caption{Schematic diagram of SMLP architecture, a typical MLP embedded with orthogonal constraints on Stiefel manifold. In the \textbf{Stiefel layer}, the output features from the previous layer are first linearly combined to generate new feature representations, processed by the orthogonal operator $\boldsymbol{\Psi}$ to strictly satisfy the orthogonal constraints on the Stiefel manifold, and output $\hat{Q}$. The \textbf{Nonnegativity Constraint} is denoted by $loss_{nonneg}:=\frac{1}{2}\| S \circ (Q \Lambda Q^T) - \text{ReLU}(S \circ (Q \Lambda Q^T)) \|_F^2$, the \textbf{Prescribed Entries Constraint} by $loss_{spec}:=\frac{1}{2}\|\Omega - \{1-\kappa_1, \kappa_1 \cdot \text{ReLU}\}(\Gamma \circ (Q \Lambda Q^T))\|_F^2$ , and the \textbf{Row Sum Constraint} by $loss_{row}:=\frac{1}{2}\|\alpha - \beta\|_F^2$. The total $\textbf{Loss} = \kappa_1 \cdot loss_{nonneg} + loss_{spec} + \kappa_2 \cdot loss_{row}$ is obtained by selecting $\kappa_1$ and $\kappa_2$. Updating parameters by minimizing the \textbf{Loss}.}\label{SMLP}
\end{figure*}
\paragraph{A unified framework for Solving SIEPs.}
Drawing from (\ref{Loss1}), (\ref{Loss2}) and (\ref{Loss3}), we propose a more unified loss function to address SIEPs under consideration.
In facts a result the above optimization problem can be written as the
following unified form:
\begin{equation*}
       \begin{aligned} 
         \min \quad &F(Q)  = \kappa_1 \cdot \frac{1}{2}\| S \circ (Q \Lambda Q^T) - \text{ReLU}(S \circ (Q \Lambda Q^T)) \|_F^2   \\ 
         & \quad \quad \quad + \frac{1}{2}\|\Omega - \{1-\kappa_1, \kappa_1 \cdot \text{ReLU}\}(\Gamma \circ (Q \Lambda Q^T))\|_F^2
         \\ & \quad \quad \quad + \kappa_2 \cdot \frac{1}{2}\|\alpha - \beta\|_F^2,\\ 
         &\text { s.t. } \quad Q \in \mathcal{O}(n),
       \end{aligned}
     \end{equation*}    
where
     $$
     \{1-\kappa_1, \kappa_1 \cdot \text{ReLU}\} = 
     \begin{cases} 
     \text{ReLU} & \text{if } \kappa_1 = 1, \\
     1 & \text{if } \kappa_1 = 0.
     \end{cases}
     $$

Therefore, for solving the SIEPs, we propose the following unified loss function :
\begin{equation}\label{Loss}
       Loss = \kappa_1 \cdot loss_{nonneg} + loss_{spec} + \kappa_2 \cdot loss_{row}
\end{equation}
where $loss_{nonneg}:=\frac{1}{2}\| S \circ (Q \Lambda Q^T) - \text{ReLU}(S \circ (Q \Lambda Q^T)) \|_F^2$
represents the constraint on the nonnegativity of matrix elements, 
$loss_{spec}:=\frac{1}{2}\|\Omega - \{1-\kappa_1, \kappa_1 \cdot \text{ReLU}\}(\Gamma \circ (Q \Lambda Q^T))\|_F^2$ 
represents the constraint on specified matrix entries, and
$loss_{row}:=\frac{1}{2}\|\alpha - \beta\|_F^2$
represents the constraint on the row sums of the matrix.
By choosing appropriate parameters $\kappa$ , the loss function can be adapted to different SIEPs:

\begin{itemize}
       \item  Case 1. $\kappa_1=0,\kappa_2=0$, the loss function (\ref{Loss}) corresponds to the PEIEP;
       \item  Case 2. $\kappa_1=1,\kappa_2=0$, the NPEIEP;
       \item  Case 3. $\kappa_1=1,\kappa_2=1$, the StIEP and GStIEP.       
\end{itemize}

Specifically, let \(\Lambda\) be the diagonal matrix with \(\lambda_1, \ldots, \lambda_n\) as its diagonal elements. 
Then Cases 1, 2 and 3 correspond to the symmetric form of the SIEPs, respectively.
\subsection{MLP with the Stiefel Layer (SMLP)}
According to (\ref{Sti}), the orthogonality constraint in (\ref{Loss}) resides within the Stiefel manifold. 
Therefore, we employ matrix orthogonal decomposition methods (QR decomposition in (\ref{QR}) or SVD in (\ref{SVD})) within the output layer to ensure the output matrix is orthogonal, meeting the orthogonality constraint. 
Consequently, we propose incorporating a Stiefel layer into the MLP, defined as SMLP, as illustrated in Figure \ref{SMLP}.
The input of SMLP is an orthogonal matrix $Q_0 \in \mathbb{R}^{n \times n}$, and the output is again an orthogonal matrix $Q \in \mathbb{R}^{n \times n}$.
$\psi$ represents the orthogonal decomposition operation. 
The SMLP, denoted as $\mathcal{N} = \{N : \mathbb{R}^{n \times n} \rightarrow \mathbb{R}^{n \times n}\}$, is
represented by :
\begin{equation*} 
  N(Q_0) = \Psi * W_{\mathcal{L}} * \Phi * W_{\mathcal{L}-1} * \cdots  * \Phi * W_0(Q_0).
\end{equation*}
The parameters in the Stiefel layer remain the weight matrices of the linear layers and are updated through backpropagation similarly to the linear layers. 
During backpropagation, the gradients will propagate through the linear layer and update the weights. 
Orthogonal decomposition is a step in the forward propagation and does not affect the backpropagation process. 
After the Stiefel layer obtains the result of the linear transformation, it performs orthogonal decomposition and returns an orthogonal matrix. 
By optimizing the loss function, the parameters of the SMLP are gradually adjusted so that the output matrix $Q$ approaches the ideal orthogonal matrix, thereby minimizing the value of the loss function.
When SMLP achieves convergence, we obtain $Q$ and the associated objective matrix expressed in (\ref{pi}) . For the process of this method, one can refer to Algorithm 1.
\begin{figure*}[h]
       \begin{center}
       \begin{minipage}{\textwidth}
       \begin{algorithm}[H]
          \caption{SMLP for solving the SIEP.}
          \begin{algorithmic}[1]
            \STATE \textbf{Input:} Initial orthogonal matrix $Q_0$.
            \STATE \textbf{Output:} Orthogonal matrix $Q$.
            \STATE \textbf{Step 1:} Give $N_{\text {epoch }}$ the maximum number of epochs, the stopping criterion $\varepsilon$, $S$, $\Omega$, and $\Lambda$.
            \STATE \textbf{Step 2:} Initialize the SMLP with a random initialization of parameters $\theta(W,b)$ and $H_0$ with the input data $Q_0$.
            \FOR{Epoch = 1 \TO $N_{epoch}$}
              \FOR{$i = 1$ \TO $\mathcal{L}-1$}
                 \STATE $Z_i = W_i \cdot H_{i-1} + b_i;$
                \STATE $H_i = \Phi*(Z_i);$
              \ENDFOR
              \STATE $\textbf{Stiefel layer:}$ $\hat{Q} = \Psi*(W_{\mathcal{L}} \cdot H_{\mathcal{L}-1} + b_{\mathcal{L}});$
              \STATE Update parameters $\theta(W,b)$  by minimizing the $Loss(\theta)$ in (\ref{Loss}).
              \IF {$Loss < \varepsilon$}
                \STATE Break.
              \ENDIF
             \ENDFOR
            \RETURN $Q$
            \STATE Obtain $\textbf{M}=\Omega + S\circ (Q\Lambda Q^T)$ as the solution of the SIEP.
          \end{algorithmic}
       \end{algorithm}
       \end{minipage}
       \end{center}
       \end{figure*}
                     
\section{Experimental results}
\paragraph{Neural network setup.}
To see the numerical performance, we implement the proposed SMLP with the loss function in the unified framework \eqref{Loss} using PyTorch 2.2.2.
Our neural network architecture is denoted by $[k_1,...,k_l]$, representing the SMLP with $l$ hidden layers, where the $l$-th layer contains $k_l$ neurons and uses either tanh or ReLU activation function. For training, we use the Adam optimizer with a learning rate of 0.01, a decay rate of 0.9, a specified tolerance of $\varepsilon=1e-8$, and maximum training epoch being equal to 10000, if not specified.

\paragraph{Experimental procedure.}
To ensure the reliability of the results, each experiment was repeated 100 times. However, for brevity and readability, this paper presents only the results of a single experiment, which includes:
the target matrix and function curves under different orthogonal decomposition strategies, various activation functions, and the number of hidden layers in neural networks.

\paragraph{Notation.}
To facilitate understanding, we define the following notations: \textbf{t} (in seconds) represents the time taken for the neural network to converge to the specified accuracy \(\varepsilon\), while \(\boldsymbol{\xi}\) denotes the ratio of the number of convergence to the total number of experiments, i.e. the success/convergence rate. 
\textbf{M} is the target matrix we aim to achieve. 
The error metric  \(\boldsymbol{er}\) is defined as $er = \frac{1}{100} \sum_{i=1}^{100} (\|\sigma(M) - \sigma\|_{\infty})_i$, where \(i\) denotes the \(i\)-th trial, representing the worst average error accuracy during the experiments. 
\(\textbf{Epoch}\) refers to the average training epoch at which the neural network converges.

\subsection{SMLP for Solving the SIEPs}
In this section, we use the proposed loss function \eqref{Loss} and SMLP to solve the inverse eigenvalue problem.
We first use SMLP to solve the IEP with specified elements. 
Then, we add nonnegativity constraints on matrix elements and structural symmetry constraints, forming the SNPEIEP. 
Additionally, by imposing further constraints on the main diagonal elements of the matrix, we further develop the problem. 
Finally, we address the complex StIEP and the generalized stochastic matrix inverse eigenvalue problem mentioned earlier, by additionally imposing constraints on the sum of elements in each row of the matrix.
\subsubsection{The Symmetric IEP with Prescribed Entries}
We consider the symmetric SIEP with prescribed entries with $n = 6$.
\begin{table*}[ht]
       \centering
       \caption{Influence of different activation functions and orthogonal decomposition techniques, as well as the number of hidden layers, on training time \textbf{t}, error metric \(\boldsymbol{er}\), convergence rate \(\boldsymbol{\xi}\), and the number of epochs to convergence \(\textbf{Epoch}\) in the proposed SMLP model for solving the symmetric SIEP with prescribed entries. The \textbf{Opt} denotes the orthogonal decomposition method chosen.}
       \label{Numerical results for Example 1} 
       \begin{threeparttable}
       \renewcommand{\arraystretch}{1.2}
       \begin{tabular}{*{7}{>{\centering\arraybackslash}p{1.8cm}}}
       \toprule
       \multirow{2}{*}{\textbf{Opt}} & \multirow{2}{*}{\textbf{$\Phi$}} & \multirow{2}{*}{\textbf{SMLP}} & \multicolumn{4}{c}{$\kappa_1 = 0, \kappa_2 = 0$} \\
       \cmidrule(lr){4-7} 
        & & & t & Epoch & er & $\xi$ \\
       \midrule  
       \multirow{2}{*}{QR}  & \multirow{2}{*}{ReLU} & [20] & 0.0799 & 218 & 3.55e-5 & 100\%  \\
        & & [20,20] & 0.0772 &  169 & 3.99e-5 & 100\%\\
        \multirow{2}{*}{QR} & \multirow{2}{*}{Tanh} & [20] & 0.0830 & 206 & 3.82e-5 & 100\%  \\
        & & [20,20] & 0.0826 &  182 & 3.69e-5 & 100\%\\
       \multirow{2}{*}{SVD} & \multirow{2}{*}{ReLU} & [20] &  0.0591 & 223 & 3.96e-5 & 100\% \\
        & & [20,20] &  0.0582 &  180 & 3.67e-5 & 100\%\\
        \multirow{2}{*}{SVD} & \multirow{2}{*}{Tanh} & [20] &  \textbf{0.0566} & 206 & 4.14e-5 & 100\% \\
        &  & [20,20] &  0.0775 &  225 & 4.03e-5 & 100\%\\
       \bottomrule
       \end{tabular}
       \end{threeparttable}
       \end{table*}
\begin{figure*}[ht]
              \centering
              \begin{subfigure}[b]{0.49\textwidth}
                  \centering
                  \includegraphics[width=\textwidth]{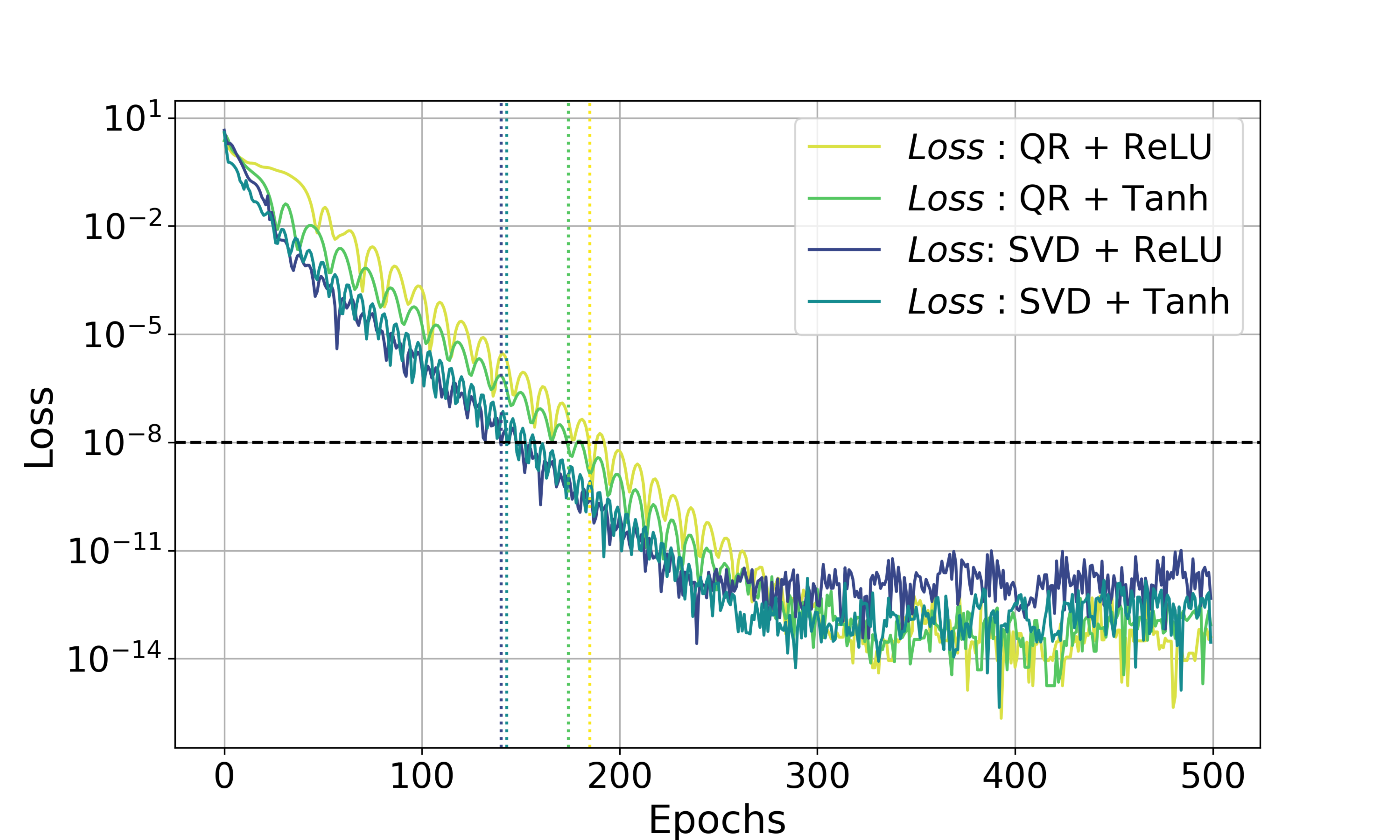}
                  \caption{Loss with one hidden layer, each with 20 neurons}
                  \label{fig:one_hidden_layer}
              \end{subfigure}
              \hfill
              \begin{subfigure}[b]{0.49\textwidth}
                  \centering
                  \includegraphics[width=\textwidth]{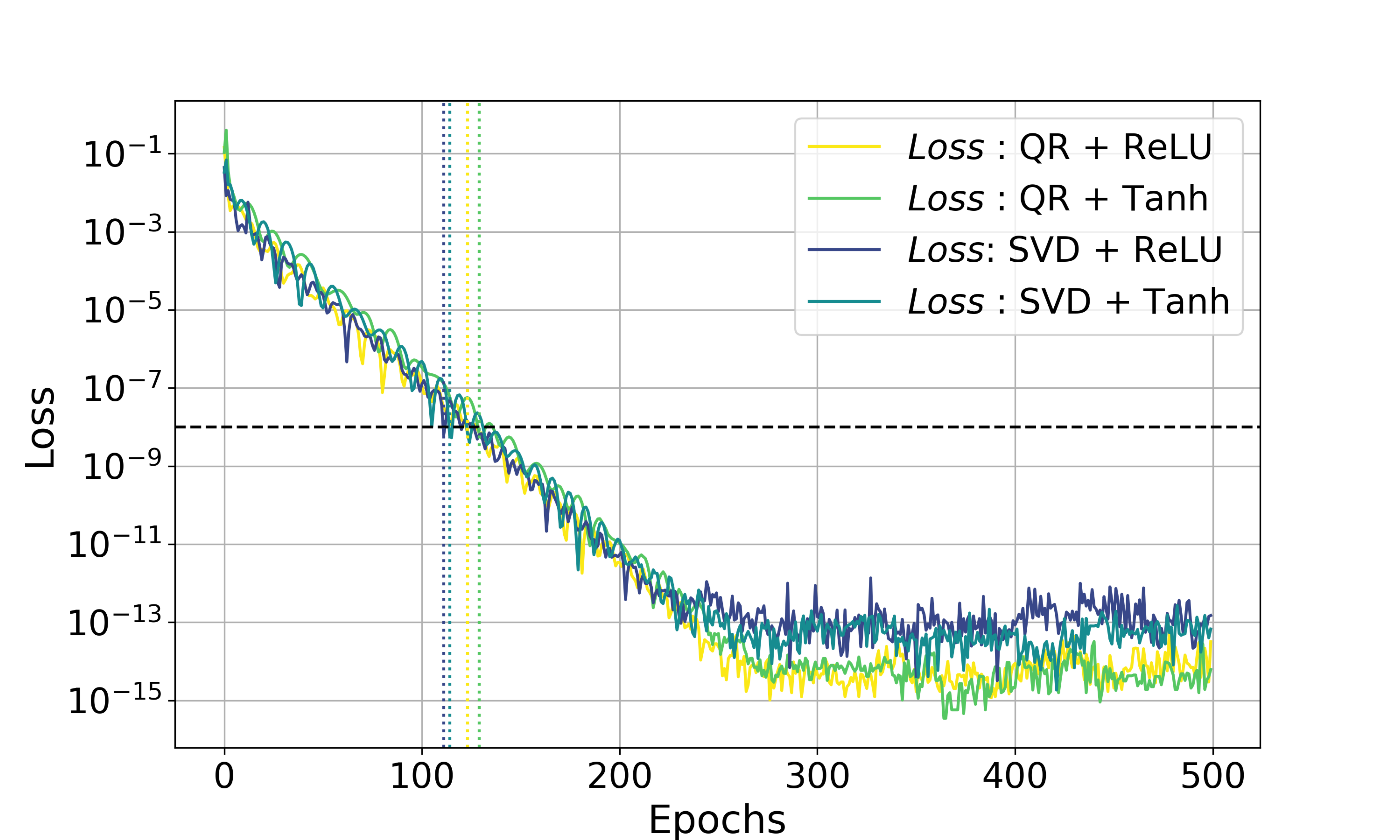}
                  \caption{Loss with two hidden layers, each with 20 neurons}
                  \label{fig:two_hidden_layers}
              \end{subfigure}
              \caption{Comparison of (a) Loss with one hidden layer; (b) Loss with two hidden layers for various combinations of different orthogonal decomposition methods and different activation functions.}
              \label{Example1}
          \end{figure*}

Let 
\(\sigma=\{2, -0.3408, 0.1046, 0.2438, -0.8483, 0.3211\}\), \(\Lambda=\text{diag}\{2, -0.3408, 0.1046, 0.2438, -0.8483, 0.3211\}\), where \(\text{diag}\) represents a diagonal matrix generated by the diagonal elements of $\sigma$, and let 
$$
       S=\left(\begin{array}{ccccccc}
              1 & \mathbf{0} & 1 &  \mathbf{0} & 1 & 1 & 1 \\
              \mathbf{0} & 1 & 1 & \mathbf{0} & 1 & 1 & 1 \\
              1 & 1 & 1 & 1 & 1 & 1 & 1 \\
              1 & 1 & 1 & 1 & 1 & 1 & 1 \\
              \mathbf{0} &  \mathbf{0} & 1 & 1 & 1 & 1 & 1 \\
              1 & 1 & 1 & 1 & 1 & 1 & 1 \\
              1 & 1 & 1 & 1 & 1 & 1 & 1
     \end{array}\right).
$$
Assume the matrix $\Omega$ as follows:
{\small
$$
\Omega=\left(\begin{array}{cccccc}
0.0000 & \textbf{0.2245} & 0.0000 & \textbf{1.3222} & 0.0000 & 0.0000 \\
\textbf{0.2245} & 0.0000 & 0.0000 & \textbf{0.4471} & 0.0000 & 0.0000 \\
0.0000 & 0.0000 & 0.0000 & 0.0000 & 0.0000 & 0.0000 \\
\textbf{1.3222} & \textbf{0.4471} & 0.0000 & 0.0000 & 0.0000 & 0.0000 \\
0.0000 & 0.0000 & 0.0000 & 0.0000 & 0.0000 & 0.0000 \\
0.0000 & 0.0000 & 0.0000 & 0.0000 & 0.0000 & 0.0000 
\end{array}\right).
$$ 
}
The numerical results are shown in Table \ref{Numerical results for Example 1}.
Additionally, the loss for different orthogonal decompositions and activation functions are presented in Figure \ref{Example1}.
The target matrix is shown as :
   {\small
   $$
   \setlength\arraycolsep{2pt} 
   \textbf{M}=\left(\begin{array}{cccccc}
       0.4732 &  \textbf{0.2245} &  0.1396 &  \textbf{1.3223} &  0.0422 & -0.1040 \\
       \textbf{0.2245} &  0.2755 & -0.0889 &  \textbf{0.4471} & -0.0271 & -0.0175 \\
       0.1396 & -0.0889 &  0.1057 &  0.1603 & -0.2188 &  0.1249 \\
       \textbf{1.3223} &  \textbf{0.4471} &  0.1603 &  0.5660 & -0.1035 & -0.0044 \\
       0.0422 & -0.0271 & -0.2188 & -0.1035 & -0.0211 &  0.2139 \\
      -0.1040 & -0.0175 &  0.1249 & -0.0044 &  0.2139 &  0.0812
   \end{array}\right).
   $$
   }
   
\textbf{Analysis.}
From Table \ref{Numerical results for Example 1}, it is evident that the training time for a configuration of one hidden layer with 20 neurons was generally shorter for the SVD compared to the QR decomposition. For instance, the SVD with ReLU  achieved a training time of 0.0591 seconds, while the QR decomposition with ReLU recorded 0.0799 seconds. This trend was consistent with both single-layer and dual-layer configurations. The convergence rate $\xi$ was 100\%, indicating that under different configurations, the SMLP could converge to a solution that satisfies the stopping criterion within the provided epochs, and it converged in fewer epochs, demonstrating the efficiency of the SMLP. 
\begin{table*}[ht]
       \centering
       \caption{Influence of different activation functions and orthogonal decomposition techniques, as well as the number of hidden layers, on training time \textbf{t}, error metric \(\boldsymbol{er}\), convergence rate \(\boldsymbol{\xi}\), and the number of epochs to convergence \(\textbf{Epoch}\) in the proposed SMLP model for solving the SNPEIEP.} 
       \label{Numerical results for Example 2} 
       \begin{threeparttable}
       \renewcommand{\arraystretch}{1.2}
       \begin{tabular}{*{7}{>{\centering\arraybackslash}p{1.8cm}}}
       \toprule
       \multirow{2}{*}{\textbf{Opt}} & \multirow{2}{*}{\textbf{$\Phi$}} & \multirow{2}{*}{\textbf{SMLP}} & \multicolumn{4}{c}{$\kappa_1$ = 1, $\kappa_2$ = 0} \\
       \cmidrule(lr){4-7} 
        & & & t & Epoch & er & $\xi$ \\
       \midrule  
       \multirow{2}{*}{QR} & \multirow{2}{*}{ReLU} & [20] & 0.0510 & 121 & 3.44e-5 & 100\%  \\
        & & [20,20] & 0.0513 & 130 & 3.63e-5 & 100\%\\
        \multirow{2}{*}{QR} & \multirow{2}{*}{Tanh} & [20] & 0.0574 & 127 & 3.39e-5 & 100\%  \\
        & & [20,20] & 0.0519 & 131 & 3.22e-5 & 100\%\\
       \multirow{2}{*}{SVD} & \multirow{2}{*}{ReLU} & [20] & 0.0382 & 120 & 3.67e-5 & 100\% \\
        & & [20,20] & \textbf{0.0339} & 128 & 3.67e-5 & 100\%\\
        \multirow{2}{*}{SVD} & \multirow{2}{*}{Tanh} & [20] & 0.0418 & 126 & 3.55e-5 & 100\% \\
        & & [20,20] & 0.0356 & 131 & 3.21e-5 & 100\%\\
       \bottomrule
       \end{tabular}
       \end{threeparttable}
       \end{table*}
       \begin{figure*}[h]
              \centering
              \begin{subfigure}[b]{0.49\textwidth}
                     \centering
                     \includegraphics[width=\textwidth]{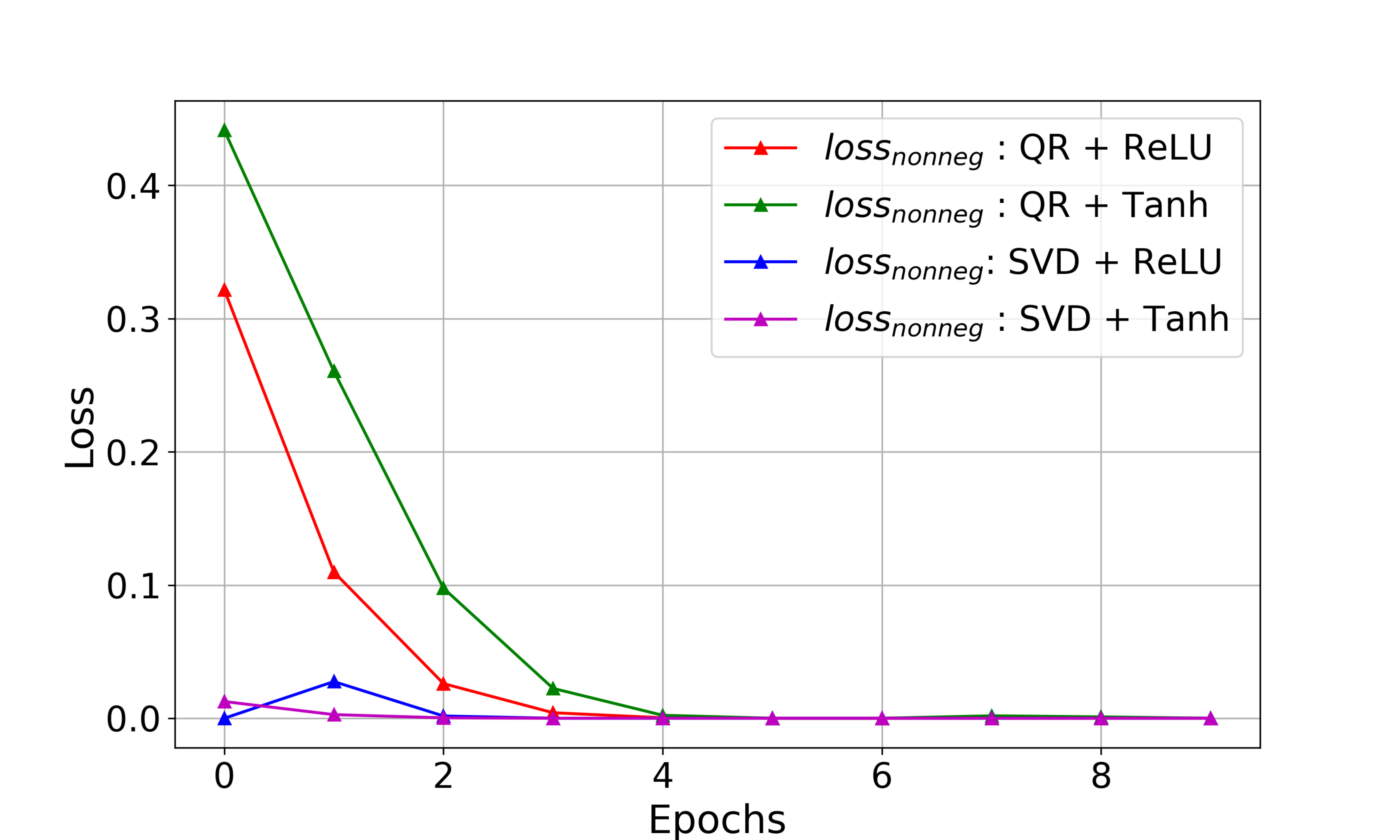}
                     \caption{$\text{Loss}_{\text{nonneg}}$ with one hidden layer, each with 20 neurons}
                 \end{subfigure}
                 \hfill
                 \begin{subfigure}[b]{0.49\textwidth}
                     \centering
                     \includegraphics[width=\textwidth]{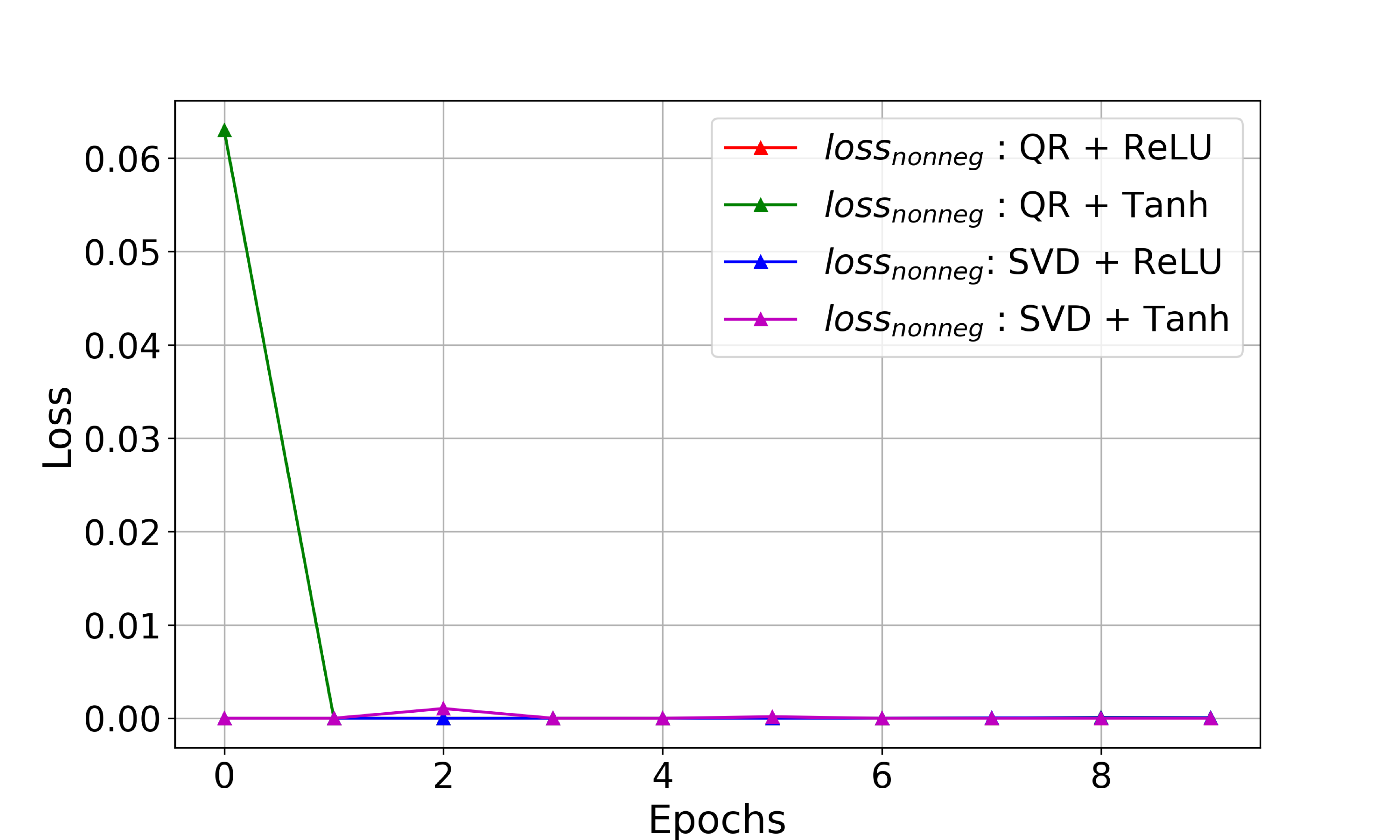}
                     \caption{$\text{Loss}_{\text{nonneg}}$ with two hidden layers, each with 20 neurons}
                 \end{subfigure}
              \caption{Comparison of (a) $\text{Loss}_{\text{nonneg}}$ with one hidden layer; (b) $\text{Loss}_{\text{nonneg}}$ with two hidden layers for various combinations of different orthogonal decomposition methods and different activation functions.}
              \label{Example2}
          \end{figure*}
From Figure \ref{Example1}, it is clear that when using SMLP to solve the PEIEP, the SVD consistently outperformed the QR decomposition in terms of training time and convergence efficiency. The choice of activation function influenced the training dynamics.
\subsubsection{The symmetric NIEP with Prescribed Entries}  
We consider the symmetric NIEP with prescribed entries with $n = 5$(SNPEIEP) \citep{chen2011isospectral}.
Given $\sigma = \{0.9568, 0.2730, 0.0253, -0.1246, -0.2352\}$, and a given subset of indices $L$ with the corresponding set of $\omega$, i.e.,
{
$$
\Omega=\left(\begin{array}{ccccc}
    \textbf{0.0596} & 0.0000 & \textbf{0.2015} & 0.0000 & 0.0000 \\ 
    0.0000 & \textbf{0.2833} & 0.0000 & \textbf{0.2116} & 0.0000 \\
    \textbf{0.2015} & 0.0000 & 0.0000 & \textbf{0.1920} & 0.0000 \\ 
    0.0000 & \textbf{0.2116} & \textbf{0.1920} & 0.0000 & 0.0000 \\ 
    0.0000 & 0.0000 & 0.0000 & 0.0000 & 0.0000
\end{array}\right).
$$
}
Let $\Lambda =\text{diag}\{0.9568,0.2730,0.0253.-0.1246.-0.2352\}$, and
{ 
$$
       S = \left(\begin{array}{cccccc}
              \textbf{0} & 1 & \textbf{0} & 1 & 1 \\
              1 & \textbf{0} & 1 & \textbf{0} & 1 \\
              \textbf{0} & 1 & 1 & \textbf{0} & 1 \\
              1 & \textbf{0} & \textbf{0} & 1 & 1 \\
              1 & 1 & 1 & 1 & 1
     \end{array}\right),
$$
}

The numerical results are shown in Table \ref{Numerical results for Example 2}.
Figure \ref{Example2} shows the loss value variation with respect to the increasing training epoch and the target matrix $\textbf{M}$ is 
$$
\textbf{M} = \left(\begin{array}{cccccc}
  \textbf{0.0596} & 0.0611 & \textbf{0.2015} & 0.1917 & 0.1746 \\
  0.0611 & \textbf{0.2833} & 0.3773 & \textbf{0.2116} & 0.0297 \\
  \textbf{0.2015} & 0.3773 & 0.1480 & \textbf{0.1920} & 0.1450 \\
  0.1917 & \textbf{0.2116} & \textbf{0.1920} & 0.3342 & 0.2493 \\
  0.1746 & 0.0297 & 0.1450 & 0.2493 & 0.0702
\end{array}\right).
$$

\textbf{Analysis.}
From Table \ref{Numerical results for Example 2}, it is evident that in both single hidden layer and dual hidden layer configurations, the SVD  generally had shorter training times than the QR decomposition. 
The convergence rate $\xi$ was consistently 100\% across all configurations, indicating successful convergence within the provided epochs for all methods. 
From Figure  \ref{Example2}, it is evident that our proposed loss function $loss_{nonneg}$ for nonnegative matrix elements has strong constraints, ensuring nonnegativity within a few epochs.
\subsubsection{The SIEP of Euclidean Distance Matrix } 
In this section, we consider the EDM inverse eigenvalue problem with $n = 7$\citep{nazari2014inverse}. 
\begin{table*}[ht]
       \centering
       \caption{Influence of different activation functions and orthogonal decomposition techniques, as well as the number of hidden layers, on training time \textbf{t}, error metric \(\boldsymbol{er}\), convergence rate \(\boldsymbol{\xi}\), and the number of epochs to convergence \(\textbf{Epoch}\) in the proposed SMLP model for solving the SIEP of Euclidean Distance Matrix.} 
       \label{Numerical results for Example 3.2} 
       \begin{threeparttable}
       \renewcommand{\arraystretch}{1.2}
       \begin{tabular}{*{7}{>{\centering\arraybackslash}p{1.8cm}}}
       \toprule
       \multirow{2}{*}{\textbf{Opt}} & \multirow{2}{*}{\textbf{$\Phi$}} & \multirow{2}{*}{\textbf{SMLP}} & \multicolumn{4}{c}{$\kappa_1$ = 1, $\kappa_2$ = 0} \\
       \cmidrule(lr){4-7} 
        & & & t & Epoch & er & $\xi$ \\
       \midrule  
       \multirow{2}{*}{QR} & \multirow{2}{*}{ReLU} & [20] & 0.0683 & 181 & 2.42e-5 & 100\%  \\
        & & [20,20] & 0.0749 & 168 & 2.66e-5 & 100\%\\
        \multirow{2}{*}{QR} & \multirow{2}{*}{Tanh} & [20] & 0.0695 & 184 & 2.37e-5 & 100\%  \\
        & & [20,20] & 0.0799 & 176 & 2.45e-5 & 100\%\\
       \multirow{2}{*}{SVD} & \multirow{2}{*}{ReLU} & [20] & 0.0510 & 181 & 2.59e-5 & 100\% \\
        & & [20,20] & 0.0550 & 167 & 3.04e-5 & 100\%\\
        \multirow{2}{*}{SVD} & \multirow{2}{*}{Tanh} & [20] & \textbf{0.0499} & 182 & 2.76e-5 & 100\% \\
        & & [20,20] & 0.0572 & 171 & 2.52e-5 & 100\%\\
       \bottomrule
       \end{tabular}
       \end{threeparttable}
       \end{table*}
       \begin{figure*}[ht]
       \centering
       \begin{subfigure}[b]{0.49\textwidth}
           \centering
           \includegraphics[width=\textwidth]{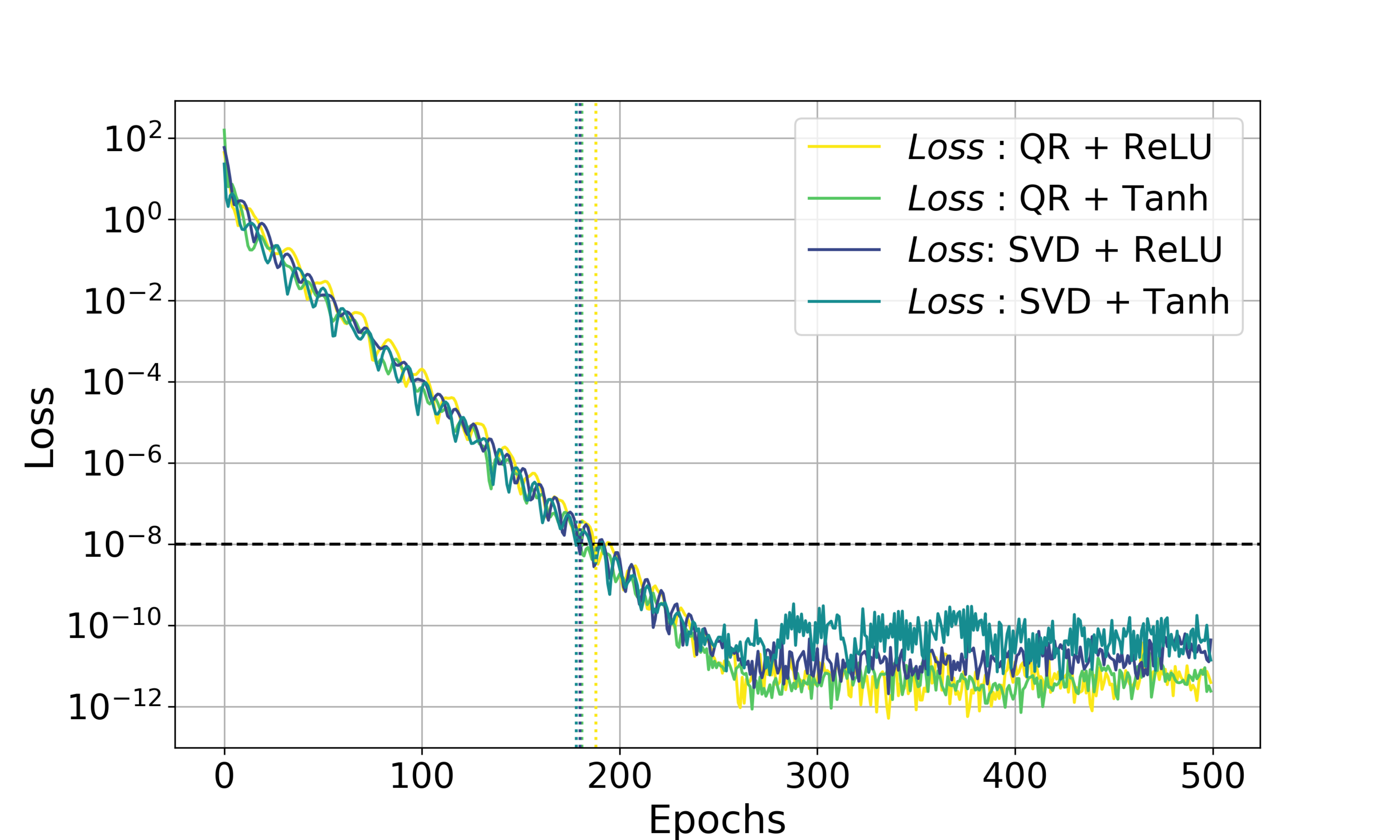}
           \caption{Loss with one hidden layer, each with 20 neurons}
       \end{subfigure}
       \hfill
       \begin{subfigure}[b]{0.49\textwidth}
           \centering
           \includegraphics[width=\textwidth]{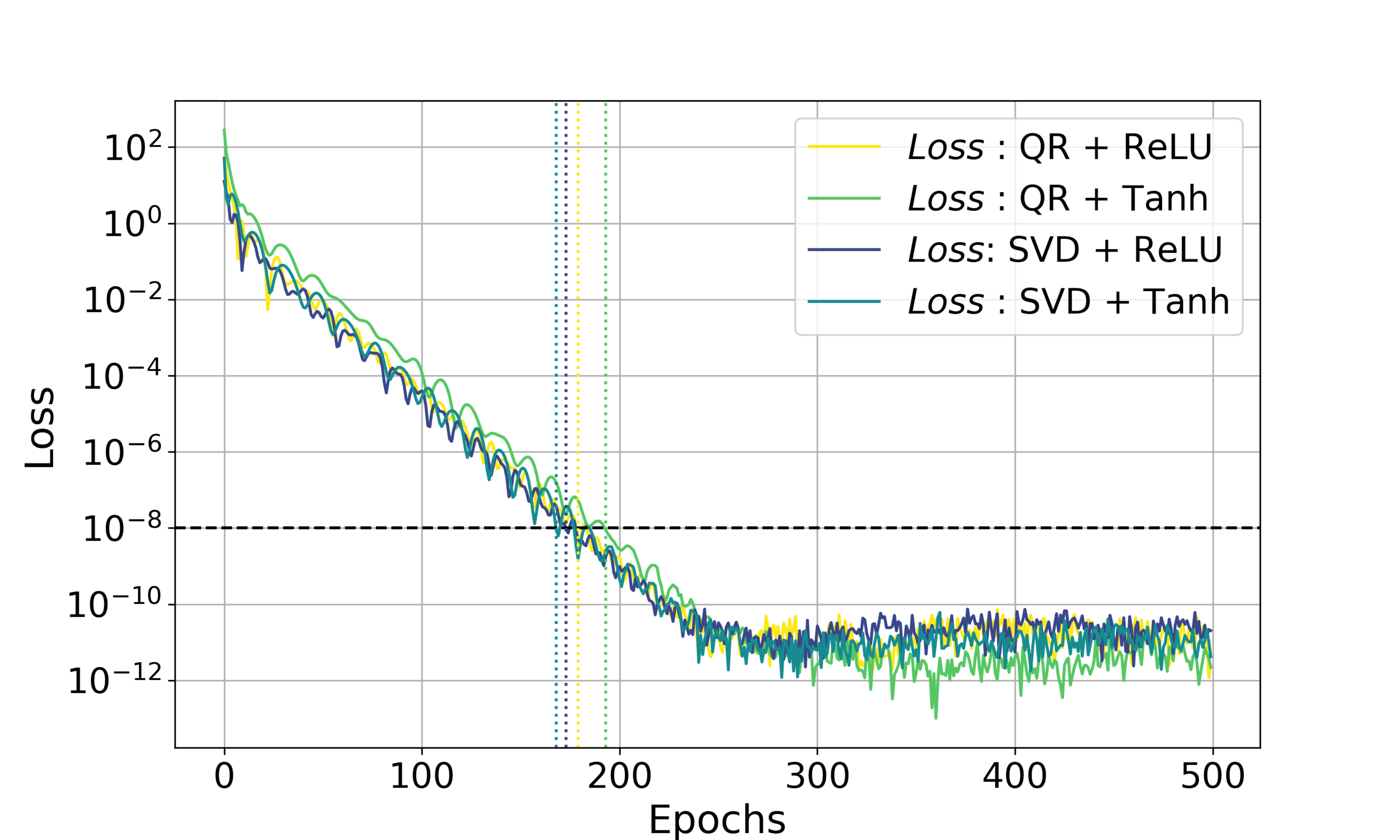}
           \caption{Loss with two hidden layers, each with 20 neurons}
       \end{subfigure}
       \caption{Comparison of (a) Loss with one hidden layer; (b) Loss with two hidden layers for various combinations of different orthogonal decomposition methods and different activation functions.}
       \label{Example3.2}
   \end{figure*}
Let $\Lambda=\text{diag}\{21,-1,-2,-3,-4,-5,-6\}$ and
$$
       S=\left(\begin{array}{ccccccc}
              \textbf{0} & 1 & 1 & 1 & 1 & 1 & 1 \\
              1 & \textbf{0} & 1 & 1 & 1 & 1 & 1 \\
              1 & 1 & \textbf{0} & 1 & 1 & 1 & 1 \\
              1 & 1 & 1 & \textbf{0} & 1 & 1 & 1 \\
              1 & 1 & 1 & 1 & \textbf{0} & 1 & 1 \\
              1 & 1 & 1 & 1 & 1 & \textbf{0} & 1 \\
              1 & 1 & 1 & 1 & 1 & 1 & \textbf{0}
     \end{array}\right).
$$
The numerical results are shown in Table \ref{Numerical results for Example 3.2}.
\begin{table*}[ht]
       \centering
       \caption{Influence of different activation functions and orthogonal decomposition techniques, as well as the number of hidden layers, on training time \textbf{t}, error metric \(\boldsymbol{er}\), convergence rate \(\boldsymbol{\xi}\), and the number of epochs to convergence \(\textbf{Epoch}\) in the proposed SMLP model for solving the StIEP.} 
       \label{Numerical results for Example 4.1} 
       \begin{threeparttable}
       \renewcommand{\arraystretch}{1.2}
       \begin{tabular}{*{7}{>{\centering\arraybackslash}p{1.8cm}}}
       \toprule
       \multirow{2}{*}{\textbf{Opt}} & \multirow{2}{*}{\textbf{$\Phi$}} & \multirow{2}{*}{\textbf{SMLP}} & \multicolumn{4}{c}{$\kappa_1$ = 1, $\kappa_2$ = 1} \\
       \cmidrule(lr){4-7} 
        & & & t & Epoch & er & $\xi$ \\
       \midrule  
       \multirow{2}{*}{QR} & \multirow{2}{*}{ReLU} & [30] & 0.3481 & 739 & 3.71e-5 & 99\%  \\
        & & [30,30] & 0.2881 & 523 & 3.54e-5 & 100\%\\
        \multirow{2}{*}{QR} & \multirow{2}{*}{Tanh} & [30] & 0.3295 & 738 & 3.80e-5 & 99\%  \\
        & & [30,30] & 0.4311 & 745 & 3.76e-5 & 100\%\\
       \multirow{2}{*}{SVD} & \multirow{2}{*}{ReLU} & [30] & \textbf{0.2688} & 815 & 3.93e-5 & 100\% \\
        & & [30,30] & 0.2860 & 682 & 3.80e-5 & 100\%\\
        \multirow{2}{*}{SVD} & \multirow{2}{*}{Tanh} & [30] & 0.2897 & 851 & 3.71e-5 & 100\% \\
        & & [30,30] & 0.3695 & 896 & 3.62e-5 & 100\%\\
       \bottomrule
       \end{tabular}
       \end{threeparttable}
       \end{table*}
\begin{figure*}[h]
              \centering
              \begin{subfigure}[b]{0.49\textwidth}
                  \centering
                  \includegraphics[width=\textwidth]{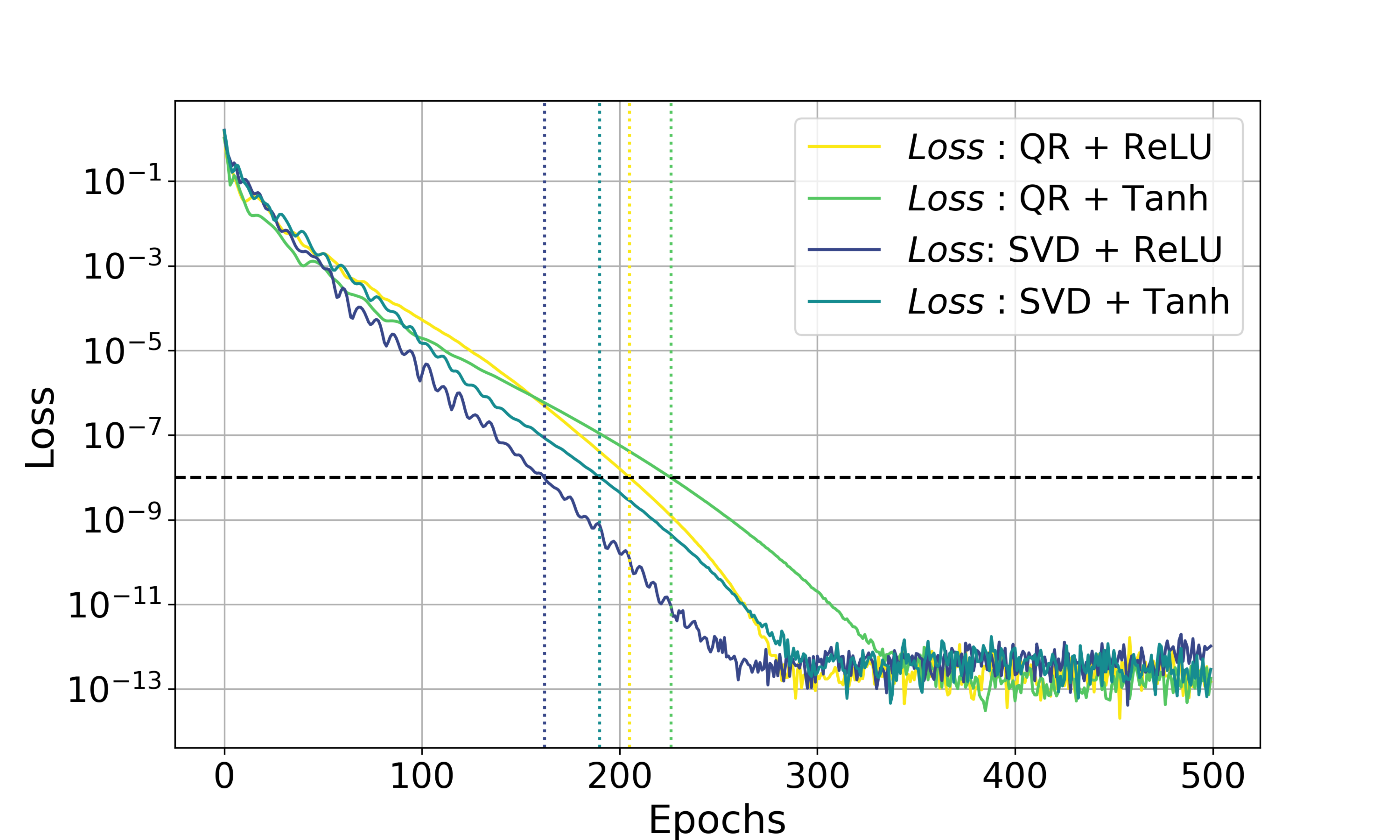}
                  \caption{Loss with one hidden layer, each with 30 neurons}
              \end{subfigure}
              \hfill
              \begin{subfigure}[b]{0.49\textwidth}
                  \centering
                  \includegraphics[width=\textwidth]{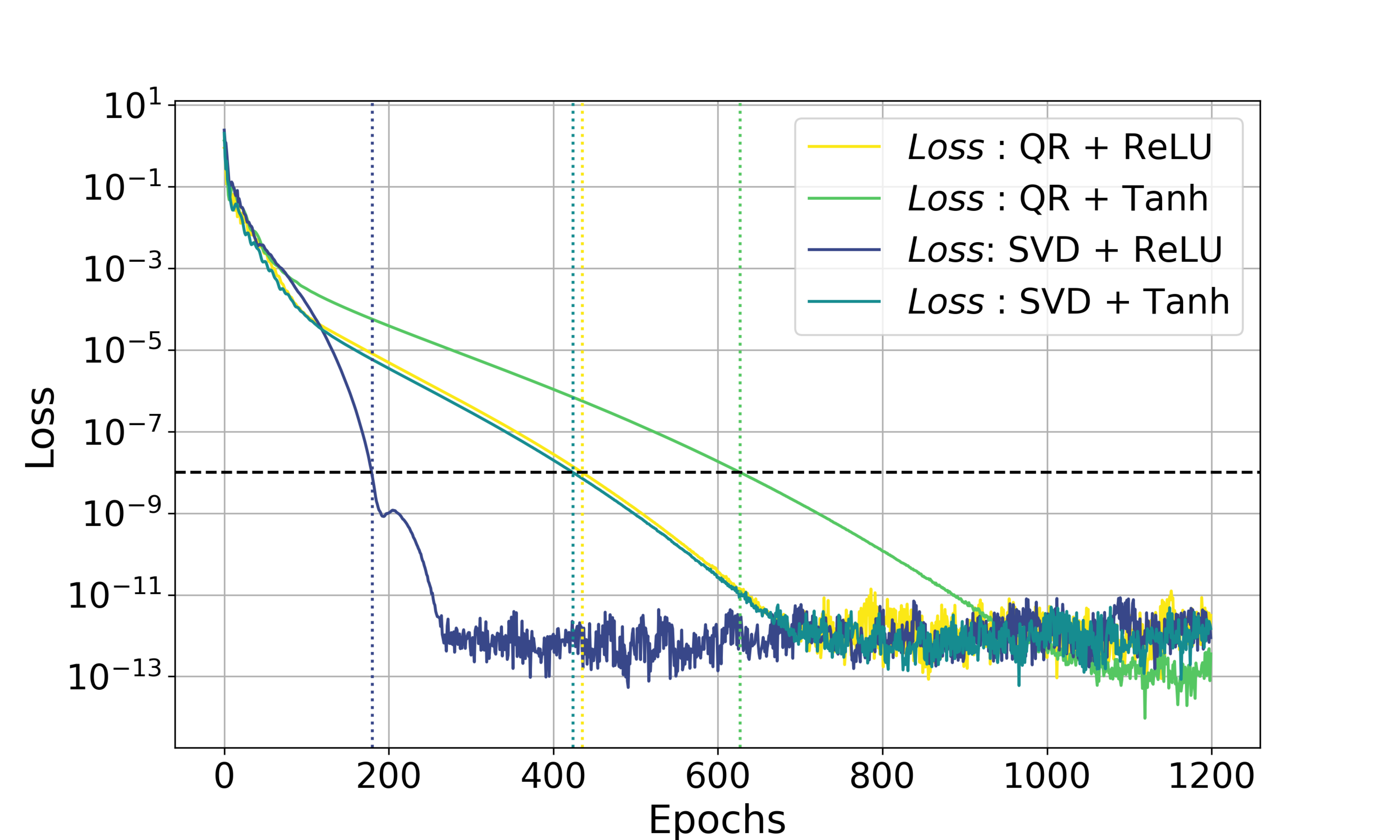}
                  \caption{Loss with two hidden layers, each with 30 neurons}
              \end{subfigure}
              \caption{Comparison of (a) Loss with one hidden layer; (b) Loss with two hidden layers for various combinations of different orthogonal decomposition methods and different activation functions.}
              \label{Example4.1}
          \end{figure*}
We selected one example, yielding the  matrix labeled as:
{\small
\setlength\arraycolsep{2pt}
$$
\textbf{M}=\left(\begin{array}{ccccccc}
       \textbf{0.0000} & 3.1681 & 4.3890 & 2.6538 & 2.3717 & 3.3239 & 4.0711 \\
       3.1681 & \textbf{0.0000} & 3.6726 & 1.8756 & 2.3498 & 4.3139 & 4.4646 \\
       4.3890 & 3.6726 & \textbf{0.0000} & 2.8574 & 3.2134 & 5.0053 & 5.6735 \\
       2.6538 & 1.8756 & 2.8574 & \textbf{0.0000} & 1.3086 & 2.7540 & 2.1672 \\
       2.3717 & 2.3498 & 3.2134 & 1.3086 & \textbf{0.0000} & 4.0498 & 3.0063 \\
       3.3239 & 4.3139 & 5.0053 & 2.7540 & 4.0498 & \textbf{0.0000} & 4.4916 \\
       4.0711 & 4.4646 & 5.6735 & 2.1672 & 3.0063 & 4.4916 & \textbf{0.0000}
   \end{array}\right).
$$
}

\textbf{Analysis.}
From Figure \ref{Example3.2} and Table  \ref{Numerical results for Example 3.2}, it can be seen that the experimental results highlight the superior performance of SVD method in terms of computational efficiency. 
Furthermore, SVD with tanh activation function demonstrated more consistent performance across different configurations. 
The findings suggest that SVD with tanh activation function constitutes a robust choice for achieving optimal numerical experiment results in SMLP models.

\subsubsection{The IEP with Row Sum Equality Constraints}  
In this example, we consider the StIEP  and the generalized stochastic matrix inverse eigenvalue problem, which requires that the sum of the elements in each row of the matrix must equal a specified value. As previously discussed, the StIEP is inherently complex.

\paragraph{The stochastic inverse eigenvalue problem.}
We consider the StIEP with $n = 5$\citep{chu1998numerical}.
Given the spectrum $\sigma = \{1.0000,-0.2608,0.5046,0.6438,-0.4483\}$.
Let $\Lambda=\text{diag}\{\linebreak[3]1.0000, -0.2608, 0.5046, 0.6438, -0.4483\}$ and 
$$
S = \left(\begin{array}{cccccc}
       1 & 1 & \textbf{0} & \textbf{0} & 1 \\
1 & 1 & 1 & \textbf{0} & \textbf{0} \\
\textbf{0} & 1 & 1 & 1 & \textbf{0} \\
\textbf{0} & \textbf{0} & 1 & 1 & 1 \\
1 & \textbf{0} & \textbf{0} & 1 & 1 
   \end{array}\right).
$$
The numerical results are shown in Table \ref{Numerical results for Example 4.1}.
And we obtained the following result:
          $$
          \textbf{M} = \left(\begin{array}{cccccc}
                 0.1916 & 0.2688 & \textbf{0.0000} & \textbf{0.0000} & 0.5396 \\
                 0.2688 & 0.2339 & 0.4972 & \textbf{0.0000} & \textbf{0.0000} \\
                 \textbf{0.0000} & 0.4972 & 0.2319 & 0.2709 & \textbf{0.0000} \\
                 \textbf{0.0000} & \textbf{0.0000} & 0.2709 & 0.5253 & 0.2038 \\
                 0.5396 & \textbf{0.0000} & \textbf{0.0000} & 0.2038 & 0.2566 \\
               \end{array}\right).
          $$

\textbf{Analysis.}
From Table \ref{Numerical results for Example 4.1}, it can be seen that in the single hidden layer configuration, SVD combined with ReLU activation achieved the lowest training time 0.2688 seconds. 
The error metric \(\boldsymbol{er}\) remains consistent across different experimental configurations, demonstrating the practical applicability of this method. Except for QR decomposition with ReLU and tanh activations in a single hidden layer configuration showing $\xi=99\%$, most methods exhibit $\xi=100\%$, a hundred percent success rate. 
SVD methods generally require more epochs to converge, especially with two hidden layers. 
From Figure \ref{Example4.1}, it can be seen that the single hidden layer configuration exhibits faster initial convergence compared to two hidden layers. 
SVD combined with ReLU consistently outperforms other methods in maintaining low loss values, particularly in the single hidden layer configuration. 
ReLU activation generally provides faster convergence and lower final loss compared to tanh. 
This trend is consistent across both QR and SVD optimizations. 
The results indicate that SVD combined with ReLU is a reliable choice for achieving optimal numerical experiment results in single-layer SMLP models with 30 neurons per layer.

\paragraph{The IEP of generalized stochastic matrix.} 
We consider the generalized stochastic matrix eigenvalue problem with $n = 7$ \citep{soto2006realizability}.
Given the spectrum $\sigma=\{8,6,3,3,-5,-5,-5,-5\}$. The aim is to construct a symmetric nonnegative matrix such that the sum of elements in each row is 8.
Let $\Lambda=\text{diag}\{8, 6, 3, 3, -5, -5, -5, -5\}$ and 
{
$$
S=\left(\begin{array}{cccccccc}
       \textbf{0} & 1 & 1 & 1 & 1 & 1 & 1 & 1 \\
       1 & \textbf{0} & 1 & 1 & 1 & 1 & 1 & 1 \\
       1 & 1 & \textbf{0} & 1 & 1 & 1 & 1 & 1 \\
       1 & 1 & 1 & \textbf{0} & 1 & 1 & 1 & 1 \\
       1 & 1 & 1 & 1 & \textbf{0} & 1 & 1 & 1 \\
       1 & 1 & 1 & 1 & 1 & \textbf{0} & 1 & 1 \\
       1 & 1 & 1 & 1 & 1 & 1 & \textbf{0} & 1 \\
       1 & 1 & 1 & 1 & 1 & 1 & 1 & \textbf{0}
   \end{array}\right).
$$
}
\begin{table*}[ht]
       \centering
       \caption{Influence of different activation functions and orthogonal decomposition techniques, as well as the number of hidden layers, on training time \textbf{t}, error metric \(\boldsymbol{er}\), convergence rate \(\boldsymbol{\xi}\), and the number of epochs to convergence \(\textbf{Epoch}\) in the proposed SMLP model for solving the GStIEP.} 
       \label{Numerical results for Example 4.2}
       \begin{threeparttable}
       \renewcommand{\arraystretch}{1.2}
       \begin{tabular}{*{7}{>{\centering\arraybackslash}p{1.8cm}}}
       \toprule
       \multirow{2}{*}{\textbf{Opt}} & \multirow{2}{*}{\textbf{$\Phi$}} & \multirow{2}{*}{\textbf{SMLP}} & \multicolumn{4}{c}{$\kappa_1$ = 1, $\kappa_2$ = 1} \\
       \cmidrule(lr){4-7} 
        & & & t & Epoch & er & $\xi$ \\
       \midrule  
       \multirow{2}{*}{QR} & \multirow{2}{*}{ReLU} & [30] & 0.3006 & 652 & 3.32e-5 & 94\%  \\
        & & [30,30] & 0.2564 & 491 & 3.06e-5 & 94\%\\
        \multirow{2}{*}{QR} & \multirow{2}{*}{Tanh} & [30] & 0.2530 & 551 & 3.33e-5 & 90\%  \\
        & & [30,30] & 0.7929 & 1483 & 2.99e-5 & 89\%\\
       \multirow{2}{*}{SVD} & \multirow{2}{*}{ReLU} & [30] & 0.1716 & 496 & 3.51e-5 & 97\% \\
        & & [30,30] & \textbf{0.1699} & 416 & 3.34e-5 & 98\%\\
        \multirow{2}{*}{SVD} & \multirow{2}{*}{Tanh} & [30] & 0.1895 & 514 & 3.54e-5 & 94\% \\
        & & [30,30] & 0.5500 & 1311 & 2.92e-5 & 93\%\\
       \bottomrule
       \end{tabular}
       \end{threeparttable}
       \end{table*}
Then the result is
       {\footnotesize
              \setlength\arraycolsep{2pt} 
       $$
               \textbf{M}=\left(\begin{array}{cccccccc}
                     \textbf{0.0000} & 0.0988 & 1.5261 & 0.2293 & 1.1586 & 4.9304 & 0.0480 & 0.0090 \\
                  0.0988 & \textbf{0.0000} & 0.0754 & 1.0774 & 0.4261 & 0.5015 & 0.8226 & 4.9982 \\
                  1.5261 & 0.0754 & \textbf{0.0000} & 0.0144 & 4.9726 & 0.8496 & 0.4920 & 0.0700 \\
                  0.2293 & 1.0774 & 0.0144 & \textbf{0.0000} & 0.0160 & 0.5391 & 4.9587 & 1.1652 \\
                  1.1586 & 0.4261 & 4.9726 & 0.0160 & \textbf{0.0000} & 0.5018 & 0.4982 & 0.4267 \\
                  4.9304 & 0.5015 & 0.8496 & 0.5391 & 0.5018 & \textbf{0.0000} & 0.2636 & 0.4140 \\
                  0.0480 & 0.8226 & 0.4920 & 4.9587 & 0.4982 & 0.2636 & \textbf{0.0000} & 0.9169 \\
                  0.0090 & 4.9982 & 0.0700 & 1.1652 & 0.4267 & 0.4140 & 0.9169 & \textbf{0.0000}
              \end{array}\right).
       $$
       }
       
The numerical results are shown in Table \ref{Numerical results for Example 4.2}. 
It can be observed that for the more complex inverse eigenvalue problem of generalized random matrices, our SMLP method still maintains a high success rate of 98\%. 
Additionally, by employing the SVD method and the ReLU activation function, the neural network achieves convergence accuracy within just 0.1699 seconds with two hidden layers.
\begin{table*}[ht]
    \centering
    \caption{Influence of different activation functions and orthogonal decomposition techniques on $\|J_{2n}-\hat{J}_{2n}\|_F$ and  $\|\lambda_{2n}-\hat{\lambda}_{2n}\|_F$ in the proposed SMLP model for solving the large scale SIEP.} 
    \label{The Large-Scale SIEP} 
    \begin{threeparttable}
    \renewcommand{\arraystretch}{1.5} 
    \begin{tabular}{>{\centering\arraybackslash}p{1.5cm} >{\centering\arraybackslash}p{1.5cm} >{\centering\arraybackslash}p{2cm} >{\centering\arraybackslash}p{2cm} >{\centering\arraybackslash}p{3.5cm} >{\centering\arraybackslash}p{3.5cm}}
    \toprule
    \multirow{2}{*}{\textbf{2n}} & \multirow{2}{*}{\textbf{SMLP}} & \multirow{2}{*}{\textbf{Opt}} & \multirow{2}{*}{\textbf{$\Phi$}} & \multicolumn{2}{c}{$\kappa_1$ = 1, $\kappa_2$ = 1} \\
    \cmidrule(lr){5-6} 
     & & & & $\scriptsize \|J_{2n}-\hat{J}_{2n}\|_F$ & $\|\lambda_{2n}-\hat{\lambda}_{2n}\|_F$ \\
    \midrule  
     \multirow{4}{*}{20} & \multirow{4}{*}{[40]} & \multirow{2}{*}{QR} & ReLU & 1.6247e-06 & 3.5738e-06 \\
     & & & Tanh& 1.3355e-06 & 5.0047e-06 \\
      & & \multirow{2}{*}{SVD} & ReLU & 4.0108e-06 & 5.3986e-06 \\
     & & & Tanh & 3.0544e-06 & 3.8078e-06 \\
     \midrule
     \multirow{4}{*}{40} & \multirow{4}{*}{[80]} & \multirow{2}{*}{QR} & ReLU & 3.0994e-06 &9.8183e-06 \\
     & & & Tanh& 3.2603e-06 & 9.3459e-06 \\
      & & \multirow{2}{*}{SVD} & ReLU & 4.0329e-06  & 8.6499e-06 \\
     & & & Tanh & 5.5870e-06 & 9.5391e-06 \\
     \midrule
     \multirow{4}{*}{60} & \multirow{4}{*}{[120]} & \multirow{2}{*}{QR} & ReLU& 3.9537e-06 &8.9414e-06  \\
     & & & Tanh&8.2261e-06 & 1.2188e-05 \\
      & & \multirow{2}{*}{SVD} & ReLU & 5.3978e-06 & 1.2243e-05 \\
     & & & Tanh & 2.4037e-05 & 4.3594e-05 \\
    \bottomrule
    \end{tabular}
    \end{threeparttable}
\end{table*}
\subsection{SIEPs of Large Scale}
In the current experiment, we aimed to demonstrate the efficiency of the proposed SMLP in solving large scale SIEPs.
We consider the following Jacobi inverse eigenvalue problem \citep{hochstadt1967some,wei2015inverse}:
given an $n \times n$ Jacobi matrix $J_n$ and a set of distinct real values $\left\{\lambda_i\right\}_{i=1}^{2 n}$, construct a $2 n \times 2 n$ Jacobi matrix $J_{2 n}$ such that $\left\{\lambda_i\right\}_{i=1}^{2 n}$ are the eigenvalues of $J_{2 n}$ and $J_n$ is the $n \times n$ leading principal submatrix of $J_{2 n}$. Particularly, for the following $\text{Jacobi}$ matrix:
\begin{equation*}
J_{2n}=\left(\begin{array}{cccccc}
2 & 1 & & & &\\
1 & 2 &1 &  & & \\
 & 1 & 2 & 1 & & \\
& &\ddots & \ddots & \ddots & \\
& & & 1 & 2 & 1 \\
& &  & & 1 & 2
\end{array}\right)_{2n\times 2n},
\end{equation*}
the eigenvalues are:
\begin{equation*}
\lambda_i=2\cos \frac{ i \pi}{2 n+1}+2, \quad i=1,2, \ldots, 2 n.
\end{equation*}
Let
\begin{equation*}
\Omega=\left(\begin{array}{cc}
J_n & \mathbf{0}\\
 \mathbf{0}& \mathbf{0}
\end{array}\right),
\end{equation*}
where $\mathbf{0}$ is an $n$ by $n$ zero matrix.
The matrix $S$ is given by:
\begin{equation*}
S=\left(\begin{array}{cc}
\mathbf{0} & \mathbf{L}\\
 \mathbf{R} & \mathbf{C}
\end{array}\right),
\end{equation*}
where $\mathbf{C}$ is an $n$ by $n$ tridiagonal matrix with all diagonal and subdiagonal elements 1, and $\mathbf{L}$ and $\mathbf{R}$ are given by:
$$
\mathbf{L}=\left(\begin{array}{cccc}
0 & 0 & \cdots & 0 \\
\vdots & \vdots & \ddots & \vdots \\
0 & 0 & \cdots & 0 \\
1 & 0 & \cdots & 0
\end{array}\right),
\mathbf{R}=\left(\begin{array}{cccc}
0 & \cdots & 0 & 1 \\
0 & \cdots & 0 & 0 \\
\vdots & \ddots & \vdots &\vdots \\
0 & \cdots & 0 & 0
\end{array}\right).
$$

We opted to train the network for a fixed number of 100000 epochs. $\hat{J}_{2n}$ is obtained through our method, and $\hat{\lambda}_{2n}$ represents its eigenvalues. The numerical results are presented in Table \ref{The Large-Scale SIEP}.

\textbf{Analysis.}
The experimental results from Table \ref{The Large-Scale SIEP} demonstrate that our proposed SMLP model performs excellently in solving large scale SIEPs. Analyzing experimental data across different dimensions reveals that the combination of QR decomposition and ReLU activation function generally yields the best performance, showcasing high solution accuracy and robustness. Our method achieves low error levels for both \(\|J_{2n} - \hat{J}_{2n}\|_F\) and \(\|\lambda_{2n} - \hat{\lambda}_{2n}\|_F\) indicators, confirming its effectiveness and practicality in solving large scale problems.

\begin{figure}[ht]
       \centering
       \includegraphics[width=0.35\textwidth]{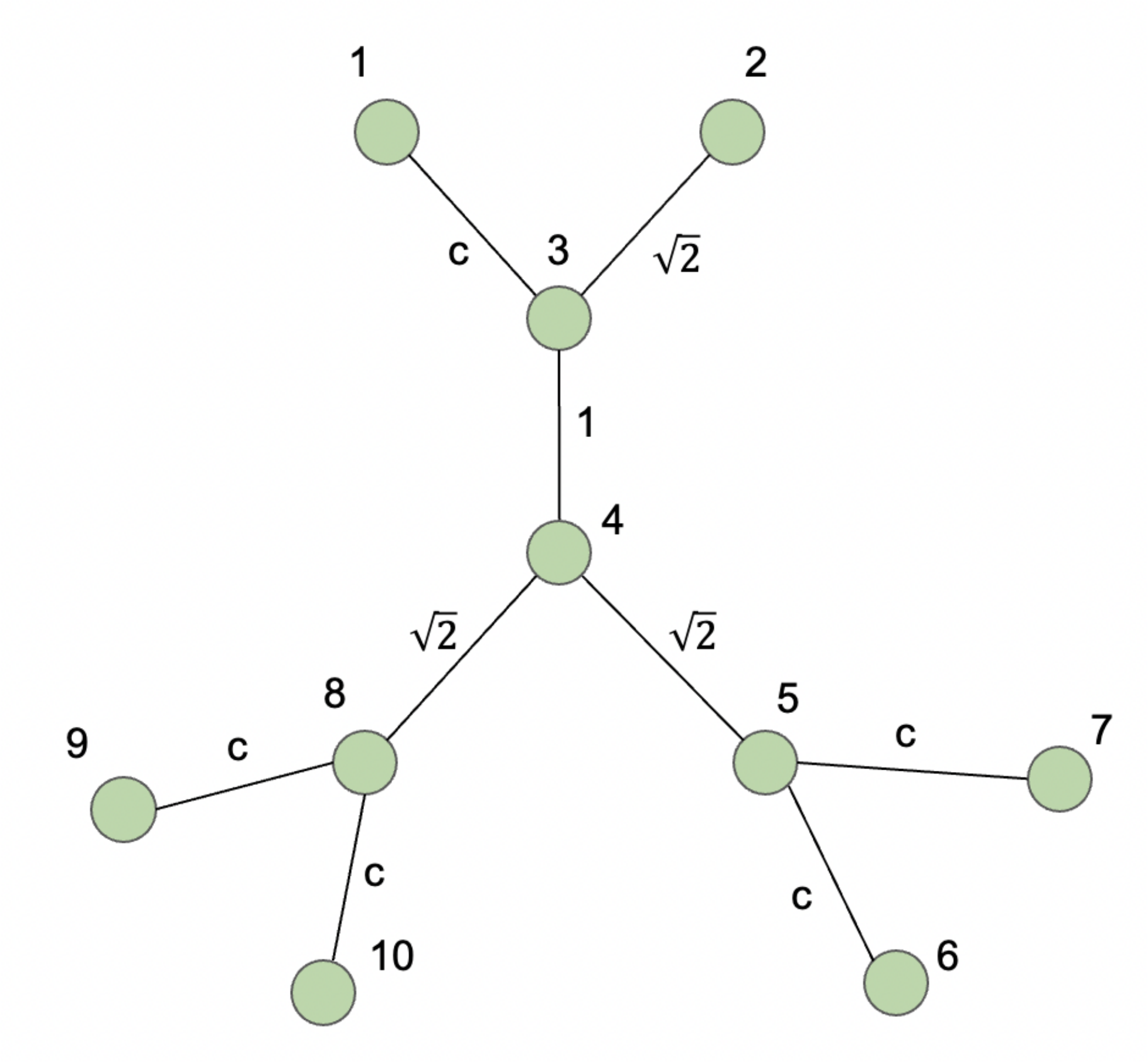} 
       \caption{Structural pattern graph of the symmetric matrix \citep{johnson2018eigenvalues,chehab2023optimization}.}\label{GPf}
\end{figure}
\subsection{Applications of SMLP in Practical Scenarios}
\paragraph{The IEP in the graph theory.}
This experiment demonstrates the application of the IEP in the field of graph theory \citep{johnson2018eigenvalues}. 
We consider a $10 \times 10$ symmetric matrix described in \citep{johnson2018eigenvalues,chehab2023optimization}, which corresponds to the undirected graph shown in Figure \ref{GPf}.
This graph is considered as the smallest possible nonlinear tree. 
The eigenvalues are set as $\sigma=\{-3, -2, -2, 0, 0, 0, 0, 2, 2, 3\}$. 
These eigenvalues and their multiplicities ensure the potential solution of the inverse eigenvalue problem. 
Combined with the constraints on the matrix elements, the corresponding matrix $\Omega$ is given as follows: 
$$
\Omega = \left(\begin{array}{cccccccccc}
       0 & 0 & \textbf{c} & 0 & 0 & 0 & 0 & 0 & 0 & 0 \\
       0 & 0 & \textbf{$\sqrt{2}$} & 0 & 0 & 0 & 0 & 0 & 0 & 0 \\
       \textbf{c} & \textbf{$\sqrt{2}$} & 0 & \textbf{1} & 0 & 0 & 0 & 0 & 0 & 0 \\
       0 & 0 & \textbf{1} & 0 & \textbf{$\sqrt{2}$} & 0 & 0 & \textbf{$\sqrt{2}$} & 0 & 0 \\
       0 & 0 & 0 & \textbf{$\sqrt{2}$} & 0 & \textbf{c} & \textbf{c} & 0 & 0 & 0 \\
       0 & 0 & 0 & 0 & \textbf{c} & 0 & 0 & 0 & 0 & 0 \\
       0 & 0 & 0 & 0 & \textbf{c} & 0 & 0 & 0 & 0 & 0 \\
       0 & 0 & 0 & \textbf{$\sqrt{2}$} & 0 & 0 & 0 & 0 & \textbf{c} & \textbf{c} \\
       0 & 0 & 0 & 0 & 0 & 0 & 0 & \textbf{c} & 0 & 0 \\
       0 & 0 & 0 & 0 & 0 & 0 & 0 & \textbf{c} & 0 & 0
     \end{array}\right),
$$
where $\textbf{c}$ indicates that the position contains a nonzero element, but the exact value is not specified in advance.
For the loss function (\ref{Loss}), we choose the parameters \(\kappa_1  = 0\) and \(\kappa _2 = 0\).
Let $\Lambda=\text{diag}\{-3, -2, -2, 0, 0, 0, 0, 2, 2, 3\}$ and 
{
$$
S=\left(\begin{array}{cccccccccc}
       \textbf{0} & \textbf{0} & 1 & \textbf{0} & \textbf{0} & \textbf{0} & \textbf{0} & \textbf{0} & \textbf{0} & \textbf{0} \\
       \textbf{0} & \textbf{0} & \textbf{0} & \textbf{0} & \textbf{0} & \textbf{0} & \textbf{0} & \textbf{0} & \textbf{0} & \textbf{0} \\
       1 & \textbf{0} & \textbf{0} & \textbf{0} & \textbf{0} & \textbf{0} & \textbf{0} & \textbf{0} & \textbf{0} & \textbf{0} \\
       \textbf{0} & \textbf{0} & \textbf{0} & \textbf{0} & \textbf{0} & \textbf{0} & \textbf{0} & \textbf{0} & \textbf{0} & \textbf{0} \\
       \textbf{0} & \textbf{0} & \textbf{0} & \textbf{0} & \textbf{0} & 1 & 1 & \textbf{0} & \textbf{0} & \textbf{0} \\
       \textbf{0} & \textbf{0} & \textbf{0} & \textbf{0} & 1 & \textbf{0} & \textbf{0} & \textbf{0} & \textbf{0} & \textbf{0} \\
       \textbf{0} & \textbf{0} & \textbf{0} & \textbf{0} & 1 & \textbf{0} & \textbf{0} & \textbf{0} & \textbf{0} & \textbf{0} \\
       \textbf{0} & \textbf{0} & \textbf{0} & \textbf{0} & \textbf{0} & \textbf{0} & \textbf{0} & \textbf{0} & 1 & 1 \\
       \textbf{0} & \textbf{0} & \textbf{0} & \textbf{0} & \textbf{0} & \textbf{0} & \textbf{0} & 1 & \textbf{0} & \textbf{0} \\
       \textbf{0} & \textbf{0} & \textbf{0} & \textbf{0} & \textbf{0} & \textbf{0} & \textbf{0} & 1 & \textbf{0} & \textbf{0}
   \end{array}\right).
$$
}
The numerical results are shown in Table \ref{Numerical results for Example Gp}.
  \begin{table*}[ht]
       \centering
       \caption{Influence of different activation functions and orthogonal decomposition techniques, as well as the number of hidden layers, on training time \textbf{t}, error metric \(\boldsymbol{er}\), convergence rate \(\boldsymbol{\xi}\), and the number of epochs to convergence \(\textbf{Epoch}\) in the proposed SMLP model for solving the IEP in the graph theory.} 
       \label{Numerical results for Example Gp} 
       \begin{threeparttable}
       \renewcommand{\arraystretch}{1.2}
       \begin{tabular}{*{7}{>{\centering\arraybackslash}p{1.8cm}}}
       \toprule
       \multirow{2}{*}{\textbf{Opt}} & \multirow{2}{*}{\textbf{$\Phi$}} & \multirow{2}{*}{\textbf{SMLP}} & \multicolumn{4}{c}{$\kappa_1$ = 0, $\kappa_2$ = 0} \\
       \cmidrule(lr){4-7} 
        & & & t & Epoch & er & $\xi$ \\
       \midrule  
       \multirow{2}{*}{QR} & \multirow{2}{*}{ReLU} & [30] & 0.2496 & 652 & 1.08e-4 & 100\%  \\
        & & [30,30] & 0.3757 & 829 & 9.73e-4 & 100\%\\
        \multirow{2}{*}{QR} & \multirow{2}{*}{Tanh} & [30] & 0.2095 & 528 & 1.18e-4 & 100\%  \\
        & & [30,30] & 0.2774 & 580 & 1.10e-4 & 100\%\\
       \multirow{2}{*}{SVD} & \multirow{2}{*}{ReLU} & [30] & \textbf{0.1529} & 556 & 1.11e-4 & 100\% \\
        & & [30,30] & 0.2112 & 571 & 1.04e-4 & 100\%\\
        \multirow{2}{*}{SVD} & \multirow{2}{*}{Tanh} & [30] & 0.1570 & 544 & 1.15e-4 & 100\% \\
        & & [30,30] & 0.2146 & 533 & 1.14e-4 & 100\%\\
       \bottomrule
       \end{tabular}
       \end{threeparttable}
       \end{table*}
       
\textbf{Analysis.}
Based on the numerical results in Table \ref{Numerical results for Example Gp}, we can conclude that the SMLP method is effective for solving the IEP in graph theory. 
This is demonstrated by the metric $\xi$ in the Table \ref{Numerical results for Example Gp}, where the value of $\xi$ reaches 100\% under different configurations. 
Additionally, the data indicate that using the combination of the SVD method and the ReLU activation function with a single hidden layer achieves convergence to the given stopping criterion faster than other combinations.
\begin{figure*}[h]
              \centering
              \begin{subfigure}[b]{0.49\textwidth}
                  \centering
                  \includegraphics[width=\textwidth]{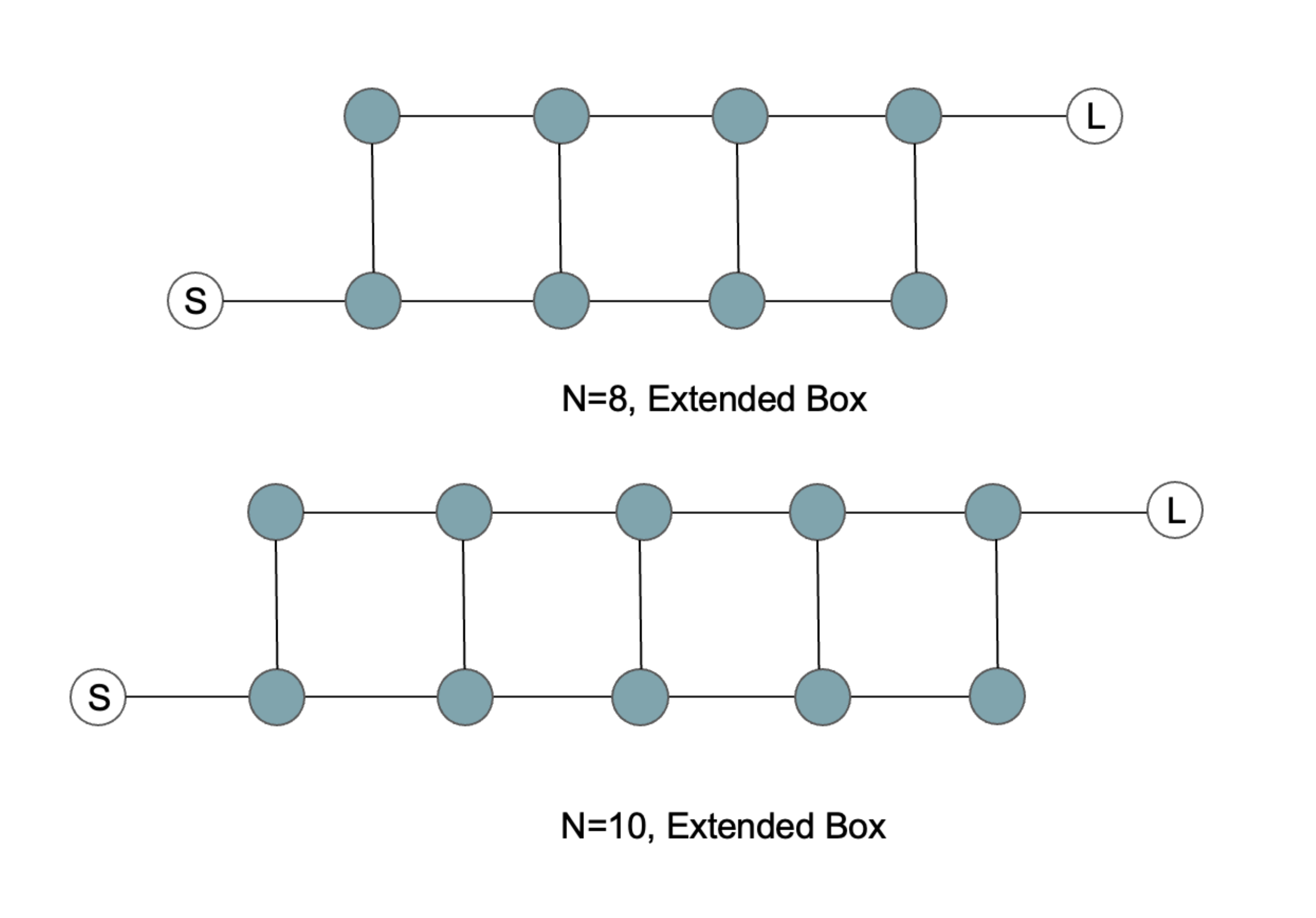}
                  \caption{Extend Box}
              \end{subfigure}
              \hfill
              \begin{subfigure}[b]{0.49\textwidth}
                  \centering
                  \includegraphics[width=\textwidth]{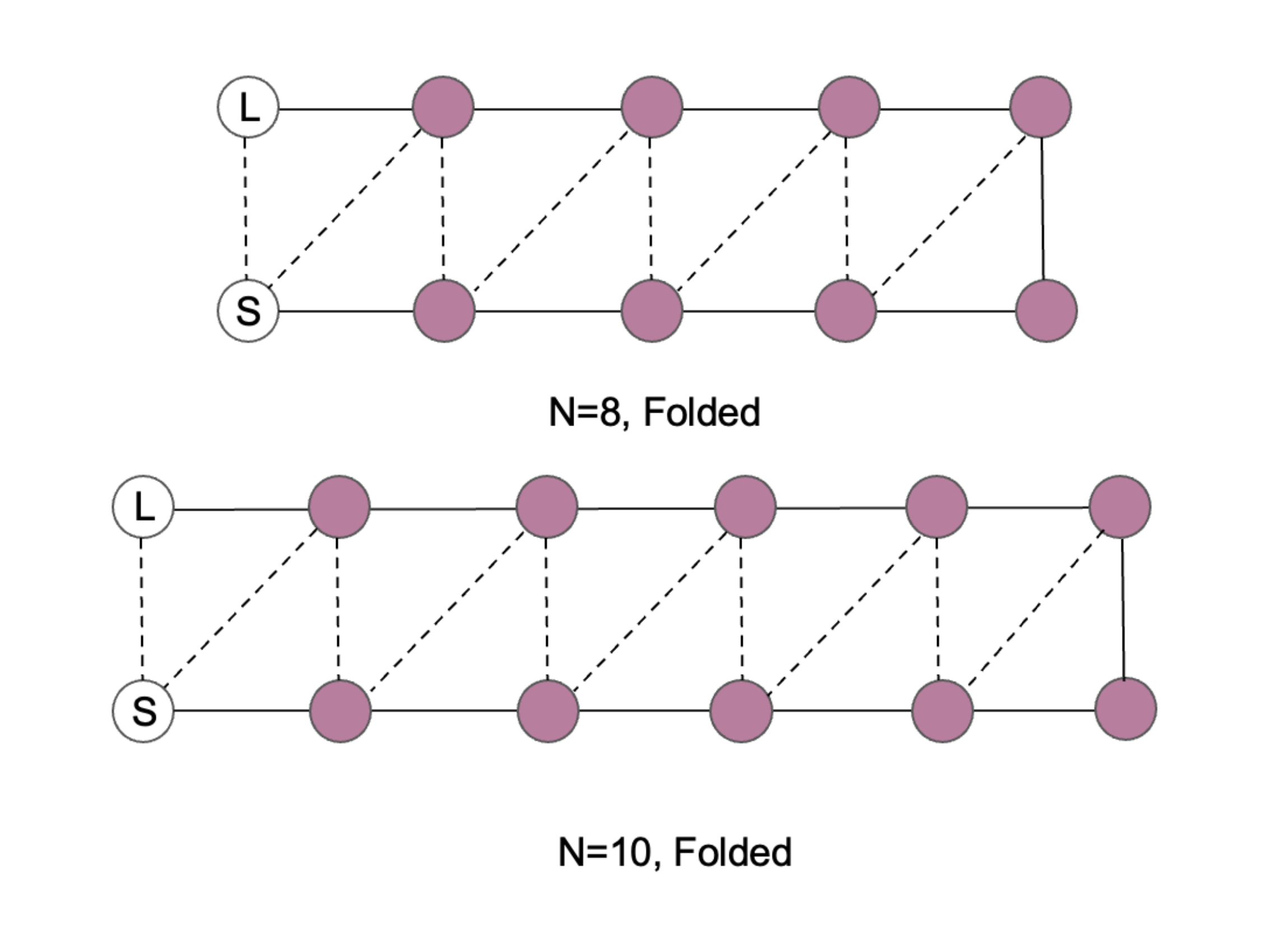}
                  \caption{Folded}
              \end{subfigure}
              \caption{Coupled Structure with Folded and Extended Box Topologies.}\label{EBFOLD}
          \end{figure*}
       \begin{figure*}[h]
              \centering
              \begin{subfigure}[b]{\textwidth}
                         \centering
       {\footnotesize 
$$
S_{\text{Folded}}=\left(\begin{array}{cccccccccc}
       \textbf{1} & \textbf{1} & 0 & 0 & 0 & 0 & 0 & 0 & 0 & \textbf{1} \\
       \textbf{1} & \textbf{1} & \textbf{1} & 0 & 0 & 0 & 0 & 0 & \textbf{1} & \textbf{1} \\
       0 & \textbf{1} & \textbf{1} & \textbf{1} & 0 & 0 & 0 & \textbf{1} & \textbf{1} & 0 \\
       0 & 0 & \textbf{1} & \textbf{1} & \textbf{1} & 0 & \textbf{1} & \textbf{1} & 0 & 0 \\
       0 & 0 & 0 & \textbf{1} & \textbf{1} & \textbf{1} & \textbf{1} & 0 & 0 & 0 \\
       0 & 0 & 0 & 0 & \textbf{1} & \textbf{1} & \textbf{1} & 0 & 0 & 0 \\
       0 & 0 & 0 & \textbf{1} & \textbf{1} & \textbf{1} & \textbf{1} & \textbf{1} & 0 & 0 \\
       0 & 0 & \textbf{1} & \textbf{1} & 0 & 0 & \textbf{1} & \textbf{1} & \textbf{1} & 0 \\
       0 & \textbf{1} & \textbf{1} & 0 & 0 & 0 & 0 & \textbf{1} & \textbf{1} & \textbf{1} \\
       \textbf{1} & \textbf{1} & 0 & 0 & 0 & 0 & 0 & 0 & \textbf{1} & \textbf{1} \\
   \end{array}\right),\quad
   S_{\text{Extended Box}}=\left(\begin{array}{cccccccccc}
       0 & \textbf{1} & 0 & 0 & 0 & 0 & 0 & 0 & 0 & 0 \\
       \textbf{1} & \textbf{1} & \textbf{1} & 0 & \textbf{1} & 0 & 0 & 0 & 0 & 0 \\
       0 & \textbf{1} & \textbf{1} & \textbf{1} & 0 & 0 & 0 & 0 & 0 & 0 \\
       0 & 0 & \textbf{1} & \textbf{1} & \textbf{1} & 0 & \textbf{1} & 0 & 0 & 0 \\
       0 & \textbf{1} & 0 & \textbf{1} & \textbf{1} & \textbf{1} & 0 & 0 & 0 & 0 \\
       0 & 0 & 0 & 0 & \textbf{1} & \textbf{1} & \textbf{1} & \textbf{1} & 0 & 0 \\
       0 & 0 & 0 & \textbf{1} & 0 & \textbf{1} & \textbf{1} & 0 & \textbf{1} & 0 \\
       0 & 0 & 0 & 0 & 0 & \textbf{1} & 0 & \textbf{1} & \textbf{1} & 0 \\
       0 & 0 & 0 & 0 & 0 & 0 & \textbf{1} & \textbf{1} & \textbf{1} & \textbf{1} \\
       0 & 0 & 0 & 0 & 0 & 0 & 0 & 0 & \textbf{1} & 0 \\
\end{array}\right).
$$
       }
       \caption{ $S$ corresponds to the Folded (left)  and Extended Box (right) topologies to a filter of order $N=8$.}
        \label{N=8}
              \end{subfigure}

                 \begin{subfigure}[b]{\textwidth}
       \centering
       {\footnotesize 
$$
      S_{\text{Folded}}=\left(\begin{array}{cccccccccccc}
              \textbf{1} & \textbf{1} & 0 & 0 & 0 & 0 & 0 & 0 & 0 & 0 & 0 & \textbf{1} \\
              \textbf{1} & \textbf{1} & \textbf{1} & 0 & 0 & 0 & 0 & 0 & 0 & 0 & \textbf{1} & \textbf{1} \\
              0 & \textbf{1} & \textbf{1} & \textbf{1} & 0 & 0 & 0 & 0 & 0 & \textbf{1} & \textbf{1} & 0 \\
              0 & 0 & \textbf{1} & \textbf{1} & \textbf{1} & 0 & 0 & 0 & \textbf{1} & \textbf{1} & 0 & 0 \\
              0 & 0 & 0 & \textbf{1} & \textbf{1} & \textbf{1} & 0 & \textbf{1} & \textbf{1} & 0 & 0 & 0 \\
              0 & 0 & 0 & 0 & \textbf{1} & \textbf{1} & \textbf{1} & \textbf{1} & 0 & 0 & 0 & 0 \\
              0 & 0 & 0 & 0 & 0 & \textbf{1} & \textbf{1} & \textbf{1} & 0 & 0 & 0 & 0 \\
              0 & 0 & 0 & 0 & \textbf{1} & \textbf{1} & \textbf{1} & \textbf{1} & \textbf{1} & 0 & 0 & 0 \\
              0 & 0 & 0 & \textbf{1} & \textbf{1} & 0 & 0 & \textbf{1} & \textbf{1} & \textbf{1} & 0 & 0 \\
              0 & 0 & \textbf{1} & \textbf{1} & 0 & 0 & 0 & 0 & \textbf{1} & \textbf{1} & \textbf{1} & 0 \\
              0 & \textbf{1} & \textbf{1} & 0 & 0 & 0 & 0 & 0 & 0 & \textbf{1} & \textbf{1} & \textbf{1} \\
              \textbf{1} & \textbf{1} & 0 & 0 & 0 & 0 & 0 & 0 & 0 & 0 & \textbf{1} & \textbf{1} \\
          \end{array}\right),\quad 
          S_{\text{Extended Box}}=\left(\begin{array}{cccccccccccc}
       0 & \textbf{1} & 0 & 0 & 0 & 0 & 0 & 0 & 0 & 0 & 0 & 0 \\
       \textbf{1} & \textbf{1} & \textbf{1} & 0 & \textbf{1} & 0 & 0 & 0 & 0 & 0 & 0 & 0 \\
       0 & \textbf{1} & \textbf{1} & \textbf{1} & 0 & 0 & 0 & 0 & 0 & 0 & 0 & 0 \\
       0 & 0 & \textbf{1} & \textbf{1} & \textbf{1} & 0 & \textbf{1} & 0 & 0 & 0 & 0 & 0 \\
    0 & \textbf{1} & 0 & \textbf{1} & \textbf{1} & \textbf{1} & 0 & 0 & 0 & 0 & 0 & 0 \\
    0 & 0 & 0 & 0 & \textbf{1} & \textbf{1} & \textbf{1} & 0 & \textbf{1} & 0 & 0 & 0 \\
0 & 0 & 0 & \textbf{1} & 0 & \textbf{1} & \textbf{1} & \textbf{1} & 0 & 0 & 0 & 0 \\
       0 & 0 & 0 & 0 & 0 & 0 & \textbf{1} & \textbf{1} & \textbf{1} & 0 & \textbf{1} & 0 \\
       0 & 0 & 0 & 0 & 0 & \textbf{1} & 0 & \textbf{1} & \textbf{1} & \textbf{1} & 0 & 0 \\
       0 & 0 & 0 & 0 & 0 & 0 & 0 & 0 & \textbf{1} & \textbf{1} & \textbf{1} & 0 \\
       0 & 0 & 0 & 0 & 0 & 0 & 0 & \textbf{1} & 0 & \textbf{1} & \textbf{1} & \textbf{1} \\
       0 & 0 & 0 & 0 & 0 & 0 & 0 & 0 & 0 & 0 & \textbf{1} & 0 \\
\end{array}\right).
$$
       }
       \caption{  $S$ corresponds to the Folded (left)  and Extended Box (right) topologies to a filter of order $N=10$.}
        \label{N=10}
                 \end{subfigure}
              \caption{The structural constraint matrices \( S \) corresponding to different topological structures for filters of different orders.}
          \end{figure*}

\paragraph{Network transformation of microwave filters.}
We consider a practical application of the inverse eigenvalue problem in engineering, focusing on the network transformation of microwave filters \citep{lamecki2004fast,kozakowski2005eigenvalue}.

As one of the basic mathematical models for studying the structure of the networks of microwave filters, the coupling matrix $\textbf{M}$ plays a crucial role in the synthesis and design process of microwave filters. 
However, the transverse coupling matrix $T$ synthesized from the transmission and reflection functions of microwave filters is not applicable in actual production \citep{young1963direct,atia1974narrow,cameron2003advanced}. 
Thus, it is necessary to obtain $\textbf{M}$ with specified topologies through matrix similarity transformations or optimization methods. 
However,  the matrix similarity transformation is limited and only applicable to certain structures. 
For a filter of order $N$, the amplitude frequency response of \( S_{11} \) and \( S_{21} \) remains constant, i.e.
\begin{equation*}
       S_{11}(T)=S_{11}(\textbf{M}),\quad S_{21}(T)=S_{21}(\textbf{M}),  
\end{equation*}
where the response functions are defined as follows: 
\begin{equation*}
       \begin{aligned}
              S_{11} =\frac{F}{E}, \quad S_{21} =\frac{P}{\epsilon E}, \quad \epsilon =\left.\frac{1}{\sqrt{10^{R L / 10}-1}} \cdot \frac{P(\omega)}{F(\omega)}\right|_{\omega= \pm 1},
       \end{aligned}
\end{equation*}
with the polynomials $E$, $F$, and $P$ being entirely determined by the order of the filter, the transmission zeros (roots of $P$), and the return loss ($RL$).

The process of deriving  \( \textbf{M} \) from  \( T \) requires that \( T \) and \( \textbf{M} \) have the same eigenvalues, which can be mathematically described as $\sigma(\textbf{M}) = \sigma(T)$. Additionally, \( \textbf{M} \) must maintain a specified structure.
Using the generalized isospectral flow method proposed by \citet{pfluger2015coupling}, we transform the solution to the aforementioned IEP  into the  optimization problem :
\begin{equation}\label{FilterLoss}
\min \quad F(Q)=\left\|Q T Q^T-S \circ\left(Q T Q^T\right)\right\|_{\mathrm{F}}^2,
\end{equation}
where the structure of $Q$ should satisfy
$$
Q=\left(\begin{array}{ccccc}
1 & 0 & \ldots & 0 & 0 \\
0 & & & & 0 \\
\vdots & & X & & \vdots \\
0 & & & & 0 \\
0 & 0 & \ldots & 0 & 1
\end{array}\right)_{(N+2)\times (N+2)}
$$
with
\begin{equation}\label{orc}
X^T X=I_{N\times N}, \quad Q^T Q=I_{(N+2)\times (N+2)}.
\end{equation}
Problem (\ref{FilterLoss}) remains an optimization problem on the Stiefel manifold, where equation (\ref{orc}) imposes a strict orthogonality constraint. 
The problem (\ref{FilterLoss}) can be reformulated as:
\begin{equation}\label{Lorc}
       \begin{aligned}
              \min \quad & F(X)=\left\|Q T Q^T-S \circ\left(Q T Q^T\right)\right\|_{\mathrm{F}}^2,\\
              \text{s.t.}  \quad &  X \in \mathcal{O}(N).
       \end{aligned}
\end{equation}
Thus, we address the problem (\ref{Lorc}) using the loss function (\ref{Loss}) with $\kappa_1=\kappa_2=0 $ and the proposed SMLP.
First, we introduce two types of topological structures as shown in Figure \ref{EBFOLD}.
The selection of matrices \( S \) are shown in Figures 10(a) and 10(b).
Based on the results of previous experiments, this study employs SVD techniques and uses ReLU as the activation function to validate the practical applications of the SMPL.
We conducted experiments with a filter of order $N$, which corresponds to an $N+2$ by $N+2$ coupling matrix.
The matrices $T$, derived from the parameters in Table \ref{Filter}, are presented in Figure \ref{TT}.
The filter parameter settings used in the experiment and the numerical results are shown in Table \ref{Filter}.
The coupling matrices \( M \) obtained through the SMLP method are shown in Figure \ref{MM8} and \ref{MM10}.
And S-parameters of the target response generated by $T$ and  $\textbf{M}$ extracted by SMLP are shown in Figure \ref{res}.
\begin{table*}[ht]
       \centering
       \caption{Filter parameter settings and numerical results for network transformation of filters. } 
       \label{Filter} 
       \begin{threeparttable}
       \renewcommand{\arraystretch}{1.3}
       \begin{tabular}{>{\centering\arraybackslash}p{1.5cm} 
                       >{\centering\arraybackslash}p{1.0cm}
                       >{\centering\arraybackslash}p{2.5cm}
                       >{\centering\arraybackslash}p{3.5cm}
                       >{\centering\arraybackslash}p{2.0cm}
                       >{\centering\arraybackslash}p{1.2cm}
                       >{\centering\arraybackslash}p{1.2cm}
                       >{\centering\arraybackslash}p{1.2cm}}
           \toprule
           \multirow{2}{*}{\textbf{Select}} & \multirow{2}{*}{\textbf{N}} & \multirow{2}{*}{\textbf{Return Loss (RL)}} & \multirow{2}{*}{\textbf{Transmission Zeroes (TZ)}} & \multirow{2}{*}{\textbf{SMLP}} & \multicolumn{3}{c}{$\kappa_1$ = 0, $\kappa_2$ = 0} \\
           \cmidrule(lr){6-8}
           & & & & & t & Epoch & $\xi$ \\
           \midrule
           \multirow{2}{*}{Folded} & 8 & 23 & [-1.3553, 1.1093, 1.2180] & [30] & 0.0871 & 297 & 100\% \\
           & 10 & 20 & [1.25, -1.8, -1.1, -1.3] & [30] & 0.1788 & 593 & 100\% \\
           \cmidrule{1-8}
           \multirow{4}{*}{Extended Box} & \multirow{2}{*}{8} & \multirow{2}{*}{23} & \multirow{2}{*}{[-1.3553, 1.1093, 1.2180]} & [30] & 0.8059 & 2973 & 100\% \\
           & & & & [30,30] & 0.8143 & 2445 & 100\% \\
           & \multirow{2}{*}{10} & \multirow{2}{*}{20} & \multirow{2}{*}{[1.25, -1.8, -1.1, -1.3]} & [30] & 1.3375 & 4736 & 75\% \\
           & & & & [30,30] & 1.5688 & 4498 & 82\% \\
           \bottomrule
       \end{tabular}
       \end{threeparttable}
   \end{table*}
\begin{figure*}[h]
       \centering
       {
       
       $$
       T=\left(\begin{array}{cccccccccc}
           0.0000 & \textbf{-0.3049} & \textbf{-0.3384} & \textbf{0.3261} & \textbf{0.4064} & \textbf{0.2741} & \textbf{-0.4140} & \textbf{-0.4059} & \textbf{0.4439} & 0.0000 \\
           \textbf{-0.3049} & \textbf{1.1742} & 0.0000 & 0.0000 & 0.0000 & 0.0000 & 0.0000 & 0.0000 & 0.0000 & \textbf{0.3049} \\
           \textbf{-0.3384} & 0.0000 & \textbf{-1.1266} & 0.0000 & 0.0000 & 0.0000 & 0.0000 & 0.0000 & 0.0000 & \textbf{0.3384} \\
           \textbf{0.3261} & 0.0000 & 0.0000 & \textbf{-1.1248} & 0.0000 & 0.0000 & 0.0000 & 0.0000 & 0.0000 & \textbf{0.3261} \\
           \textbf{0.4064} & 0.0000 & 0.0000 & 0.0000 & \textbf{1.1098} & 0.0000 & 0.0000 & 0.0000 & 0.0000 & \textbf{0.4064} \\
           \textbf{0.2741} & 0.0000 & 0.0000 & 0.0000 & 0.0000 & \textbf{-0.9343} & 0.0000 & 0.0000 & 0.0000 & \textbf{0.2741} \\
           \textbf{-0.4140} & 0.0000 & 0.0000 & 0.0000 & 0.0000 & 0.0000 & \textbf{0.6831} & 0.0000 & 0.0000 & \textbf{0.4140} \\
           \textbf{-0.4059} & 0.0000 & 0.0000 & 0.0000 & 0.0000 & 0.0000 & 0.0000 & \textbf{-0.5550} & 0.0000 & \textbf{0.4059} \\
           \textbf{0.4439} & 0.0000 & 0.0000 & 0.0000 & 0.0000 & 0.0000 & 0.0000 & 0.0000 & \textbf{0.0623} & \textbf{0.4439} \\
           0.0000 & \textbf{0.3049} & \textbf{0.3384} & \textbf{0.3261} & \textbf{0.4064} & \textbf{0.2741} & \textbf{0.4140} & \textbf{0.4059} & \textbf{0.4439} & 0.0000 \\
       \end{array}\right).
       $$
       $$
       T=\left(\begin{array}{cccccccccccc}
           0.0000 & \textbf{-0.2684} & \textbf{0.2828} & \textbf{-0.2980} & \textbf{0.3387} & \textbf{0.2042} & \textbf{-0.3234} & \textbf{-0.2899} & \textbf{0.3399} & \textbf{0.3613} & \textbf{-0.3647} & 0.0000 \\
           \textbf{-0.2684} & \textbf{-1.1113} & 0.0000 & 0.0000 & 0.0000 & 0.0000 & 0.0000 & 0.0000 & 0.0000 & 0.0000 & 0.0000 & \textbf{0.2684} \\
           \textbf{0.2828} & 0.0000 & \textbf{1.0864} & 0.0000 & 0.0000 & 0.0000 & 0.0000 & 0.0000 & 0.0000 & 0.0000 & 0.0000 & \textbf{0.2828} \\
           \textbf{-0.2980} & 0.0000 & 0.0000 & \textbf{1.0853} & 0.0000 & 0.0000 & 0.0000 & 0.0000 & 0.0000 & 0.0000 & 0.0000 & \textbf{0.2980} \\
           \textbf{0.3387} & 0.0000 & 0.0000 & 0.0000 & \textbf{-1.0780} & 0.0000 & 0.0000 & 0.0000 & 0.0000 & 0.0000 & 0.0000 & \textbf{0.3387} \\
           \textbf{0.2042} & 0.0000 & 0.0000 & 0.0000 & 0.0000 & \textbf{0.9595} & 0.0000 & 0.0000 & 0.0000 & 0.0000 & 0.0000 & \textbf{0.2042} \\
           \textbf{-0.3234} & 0.0000 & 0.0000 & 0.0000 & 0.0000 & 0.0000 & \textbf{-0.8245} & 0.0000 & 0.0000 & 0.0000 & 0.0000 & \textbf{0.3234} \\
           \textbf{-0.2899} & 0.0000 & 0.0000 & 0.0000 & 0.0000 & 0.0000 & 0.0000 & \textbf{0.7726} & 0.0000 & 0.0000 & 0.0000 & \textbf{0.2899} \\
           \textbf{0.3399} & 0.0000 & 0.0000 & 0.0000 & 0.0000 & 0.0000 & 0.0000 & 0.0000 & \textbf{0.4488} & 0.0000 & 0.0000 & \textbf{0.3399} \\
           \textbf{0.3613} & 0.0000 & 0.0000 & 0.0000 & 0.0000 & 0.0000 & 0.0000 & 0.0000 & 0.0000 & \textbf{-0.4418} & 0.0000 & \textbf{0.3613} \\
           \textbf{-0.3647} & 0.0000 & 0.0000 & 0.0000 & 0.0000 & 0.0000 & 0.0000 & 0.0000 & 0.0000 & 0.0000 & \textbf{0.0174} & \textbf{0.3647} \\
           0.0000 & \textbf{0.2684} & \textbf{0.2828} & \textbf{0.2980} & \textbf{0.3387} & \textbf{0.2042} & \textbf{0.3234} & \textbf{0.2899} & \textbf{0.3399} & \textbf{0.3613} & \textbf{0.3647} & 0.0000 \\
       \end{array}\right).
 $$
       }
       \caption{ Transverse coupling matrices $T$  corresponding to a filter of order $N=8$ (top) and $N=10$ (bottom).}
       \label{TT}
\end{figure*}
\begin{figure*}[h]
       \centering
       {
$$
        \left(\begin{array}{cccccccccc}
              \textbf{0.0000} & \textbf{1.0428} & 0.0000 & 0.0000 & 0.0000 & 0.0000 & 0.0000 & 0.0000 & 0.0000 & \textbf{0.0000} \\
              \textbf{1.0428} & \textbf{0.0106} & \textbf{-0.8623} & 0.0000 & 0.0000 & 0.0000 & 0.0000 & 0.0000 & \textbf{0.0000} & \textbf{-0.0001} \\
              0.0000 & \textbf{-0.8623} & \textbf{0.0115} & \textbf{-0.5994} & 0.0000 & 0.0000 & 0.0000 & \textbf{0.0000} & \textbf{-0.0001} & 0.0000 \\
              0.0000 & 0.0000 & \textbf{-0.5994} & \textbf{0.0133} & \textbf{-0.5356} & 0.0000 & \textbf{0.0457} & \textbf{0.1316} & 0.0000 & 0.0000 \\
              0.0000 & 0.0000 & 0.0000 & \textbf{-0.5356} & \textbf{0.0898} & \textbf{-0.3361} & \textbf{0.5673} & 0.0000 & 0.0000 & 0.0000 \\
              0.0000 & 0.0000 & 0.0000 & 0.0000 & \textbf{-0.3361} & \textbf{-0.8513} & \textbf{-0.3191} & 0.0000 & 0.0000 & 0.0000 \\
              0.0000 & 0.0000 & 0.0000 & \textbf{0.0457} & \textbf{0.5673} & \textbf{-0.3191} & \textbf{-0.0073} & \textbf{0.5848} & 0.0000 & 0.0000 \\
              0.0000 & 0.0000 & \textbf{0.0000} & \textbf{0.1316} & 0.0000 & 0.0000 & \textbf{0.5848} & \textbf{0.0114} & \textbf{-0.8623} & 0.0000 \\
              0.0000 & \textbf{0.0000} & \textbf{-0.0001} & 0.0000 & 0.0000 & 0.0000 & 0.0000 & \textbf{-0.8623} & \textbf{0.0107} & \textbf{-1.0428} \\
              \textbf{0.0000} & \textbf{-0.0001} & 0.0000 & 0.0000 & 0.0000 & 0.0000 & 0.0000 & 0.0000 & \textbf{-1.0428} & \textbf{0.0000} \\
       \end{array}\right)
$$

$$
    \left(\begin{array}{cccccccccc}
           0.0000 & \textbf{1.0428} & 0.0000 & 0.0000 & 0.0000 & 0.0000 & 0.0000 & 0.0000 & 0.0000 & 0.0000 \\
           \textbf{1.0428} & \textbf{0.0107} & \textbf{0.8530} & 0.0000 & \textbf{0.1261} & 0.0000 & 0.0000 & 0.0000 & 0.0000 & 0.0000 \\
           0.0000 & \textbf{0.8530} & \textbf{0.0112} & \textbf{-0.6783} & 0.0000 & 0.0000 & 0.0000 & 0.0000 & 0.0000 & 0.0000 \\
           0.0000 & 0.0000 & \textbf{-0.6783} & \textbf{0.0152} & \textbf{0.5334} & 0.0000 & \textbf{0.0662} & 0.0000 & 0.0000 & 0.0000 \\
           0.0000 & \textbf{0.1261} & 0.0000 & \textbf{0.5334} & \textbf{0.0226} & \textbf{-0.5983} & 0.0000 & 0.0000 & 0.0000 & 0.0000 \\
           0.0000 & 0.0000 & 0.0000 & 0.0000 & \textbf{-0.5983} & \textbf{0.0605} & \textbf{0.5475} & \textbf{0.1081} & 0.0000 & 0.0000 \\
           0.0000 & 0.0000 & 0.0000 & \textbf{0.0662} & 0.0000 & \textbf{0.5475} & \textbf{0.1384} & 0.0000 & \textbf{0.8119} & 0.0000 \\
           0.0000 & 0.0000 & 0.0000 & 0.0000 & 0.0000 & \textbf{0.1081} & 0.0000 & \textbf{-0.9805} & \textbf{-0.2905} & 0.0000 \\
           0.0000 & 0.0000 & 0.0000 & 0.0000 & 0.0000 & 0.0000 & \textbf{0.8119} & \textbf{-0.2905} & \textbf{0.0106} & \textbf{-1.0428} \\
           0.0000 & 0.0000 & 0.0000 & 0.0000 & 0.0000 & 0.0000 & 0.0000 & 0.0000 & \textbf{-1.0428} & 0.0000 \\
       \end{array}\right)
$$
       }
       \caption{The matrices $M$ corresponds to the Folded (top) and Extended Box (bottom) topologies to a filter of order $N=8$.}
       \label{MM8}
\end{figure*}
\begin{figure*}[h]
       \centering
       {\footnotesize 
$$
              \left(\begin{array}{cccccccccccc}
                     \textbf{0.0000} & \textbf{0.9823} & 0.0000 & 0.0000 & 0.0000 & 0.0000 & 0.0000 & 0.0000 & 0.0000 & 0.0000 & 0.0000 & \textbf{0.0000} \\
                         \textbf{0.9823} & \textbf{-0.0054} & \textbf{0.8078} & 0.0000 & 0.0000 & 0.0000 & 0.0000 & 0.0000 & 0.0000 & 0.0000 & \textbf{0.0000} & \textbf{0.0000} \\
                         0.0000 & \textbf{0.8078} & \textbf{-0.0056} & \textbf{-0.5798} & 0.0000 & 0.0000 & 0.0000 & 0.0000 & 0.0000 & \textbf{0.0000} & \textbf{0.0000} & 0.0000 \\
                         0.0000 & 0.0000 & \textbf{-0.5798} & \textbf{-0.0062} & \textbf{-0.5388} & 0.0000 & 0.0000 & 0.0000 & \textbf{-0.0336} & \textbf{0.0000} & 0.0000 & 0.0000 \\
                         0.0000 & 0.0000 & 0.0000 & \textbf{-0.5388} & \textbf{0.0082} & \textbf{0.5087} & 0.0000 & \textbf{-0.0852} & \textbf{-0.1246} & 0.0000 & 0.0000 & 0.0000 \\
                         0.0000 & 0.0000 & 0.0000 & 0.0000 & \textbf{0.5087} & \textbf{0.0858} & \textbf{0.2335} & \textbf{0.5727} & 0.0000 & 0.0000 & 0.0000 & 0.0000 \\
                         0.0000 & 0.0000 & 0.0000 & 0.0000 & 0.0000 & \textbf{0.2335} & \textbf{0.8871} & \textbf{-0.2632} & 0.0000 & 0.0000 & 0.0000 & 0.0000 \\
                         0.0000 & 0.0000 & 0.0000 & 0.0000 & \textbf{-0.0852} & \textbf{0.5727} & \textbf{-0.2632} & \textbf{-0.0323} & \textbf{0.5242} & 0.0000 & 0.0000 & 0.0000 \\
                         0.0000 & 0.0000 & 0.0000 & \textbf{-0.0336} & \textbf{-0.1246} & 0.0000 & 0.0000 & \textbf{0.5242} & \textbf{-0.0062} & \textbf{-0.5798} & 0.0000 & 0.0000 \\
                         0.0000 & 0.0000 & \textbf{0.0000} & \textbf{0.0000} & 0.0000 & 0.0000 & 0.0000 & 0.0000 & \textbf{-0.5798} & \textbf{-0.0056} & \textbf{-0.8078} & 0.0000 \\
                         0.0000 & \textbf{0.0000} & \textbf{0.0000} & 0.0000 & 0.0000 & 0.0000 & 0.0000 & 0.0000 & 0.0000 & \textbf{-0.8078} & \textbf{-0.0054} & \textbf{-0.9823} \\
                         \textbf{0.0000} & \textbf{0.0000} & 0.0000 & 0.0000 & 0.0000 & 0.0000 & 0.0000 & 0.0000 & 0.0000 & 0.0000 & \textbf{-0.9823} & \textbf{0.0000} \\
              \end{array}\right)
$$

              $$
                 \left(\begin{array}{cccccccccccc}
                         0.0000 & \textbf{0.9823} & 0.0000 & 0.0000 & 0.0000 & 0.0000 & 0.0000 & 0.0000 & 0.0000 & 0.0000 & 0.0000 & 0.0000 \\
                         \textbf{0.9823} & \textbf{-0.0054} & \textbf{0.7590} & 0.0000 & \textbf{0.2765} & 0.0000 & 0.0000 & 0.0000 & 0.0000 & 0.0000 & 0.0000 & 0.0000 \\
                         0.0000 & \textbf{0.7590} & \textbf{-0.0048} & \textbf{0.7579} & 0.0000 & 0.0000 & 0.0000 & 0.0000 & 0.0000 & 0.0000 & 0.0000 & 0.0000 \\
                         0.0000 & 0.0000 & \textbf{0.7579} & \textbf{-0.0094} & \textbf{-0.4685} & 0.0000 & \textbf{0.0016} & 0.0000 & 0.0000 & 0.0000 & 0.0000 & 0.0000 \\
                         0.0000 & \textbf{0.2765} & 0.0000 & \textbf{-0.4685} & \textbf{-0.0118} & \textbf{0.5202} & 0.0000 & 0.0000 & 0.0000 & 0.0000 & 0.0000 & 0.0000 \\
                         0.0000 & 0.0000 & 0.0000 & 0.0000 & \textbf{0.5202} & \textbf{-0.0181} & \textbf{-0.3375} & 0.0000 & \textbf{-0.3989} & 0.0000 & 0.0000 & 0.0000 \\
                         0.0000 & 0.0000 & 0.0000 & \textbf{0.0016} & 0.0000 & \textbf{-0.3375} & \textbf{0.5056} & \textbf{-0.3558} & 0.0000 & 0.0000 & 0.0000 & 0.0000 \\
                         0.0000 & 0.0000 & 0.0000 & 0.0000 & 0.0000 & 0.0000 & \textbf{-0.3558} & \textbf{-0.0094} & \textbf{0.4332} & 0.0000 & \textbf{0.8061} & 0.0000 \\
                         0.0000 & 0.0000 & 0.0000 & 0.0000 & 0.0000 & \textbf{-0.3989} & 0.0000 & \textbf{0.4332} & \textbf{-0.4084} & \textbf{0.3431} & 0.0000 & 0.0000 \\
                         0.0000 & 0.0000 & 0.0000 & 0.0000 & 0.0000 & 0.0000 & 0.0000 & 0.0000 & \textbf{0.3431} & \textbf{0.8815} & \textbf{0.0528} & 0.0000 \\
                         0.0000 & 0.0000 & 0.0000 & 0.0000 & 0.0000 & 0.0000 & 0.0000 & \textbf{0.8061} & 0.0000 & \textbf{0.0528} & \textbf{-0.0054} & \textbf{0.9823} \\
                         0.0000 & 0.0000 & 0.0000 & 0.0000 & 0.0000 & 0.0000 & 0.0000 & 0.0000 & 0.0000 & 0.0000 & \textbf{0.9823} & 0.0000 \\
                     \end{array}\right)
              $$
       }
       \caption{The  matrices $M$ corresponds to the Folded (top) and Extended Box (bottom) topologies to a filter of order $N=10$.}
       \label{MM10}
\end{figure*}
\begin{figure*}[h!]
       \centering
       \begin{subfigure}[b]{0.24\linewidth}
           \centering
           \includegraphics[width=\linewidth]{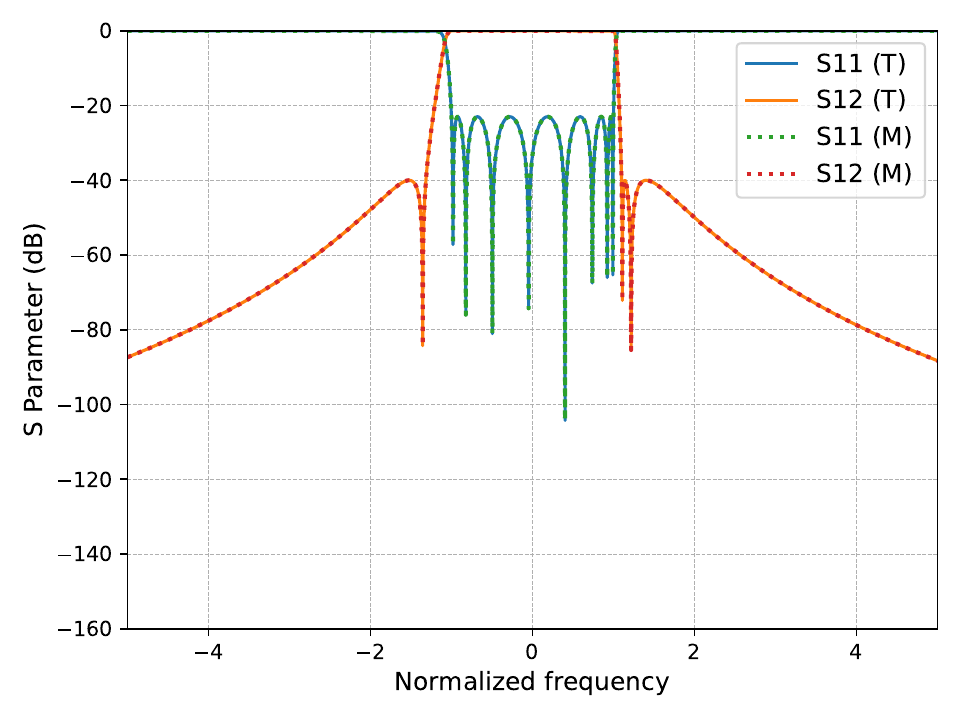}
           \caption{N=8,Folded}
       \end{subfigure}
       \begin{subfigure}[b]{0.24\linewidth}
           \centering
           \includegraphics[width=\linewidth]{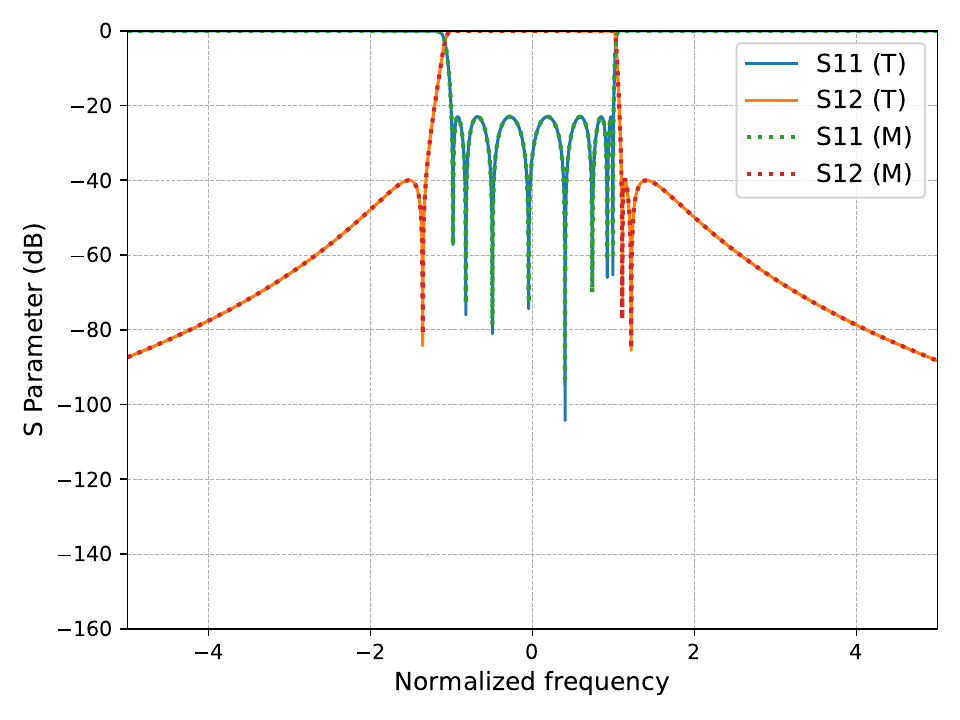}
           \caption{N=8,Extended Box}
       \end{subfigure}
       \begin{subfigure}[b]{0.24\linewidth}
           \centering
           \includegraphics[width=\linewidth]{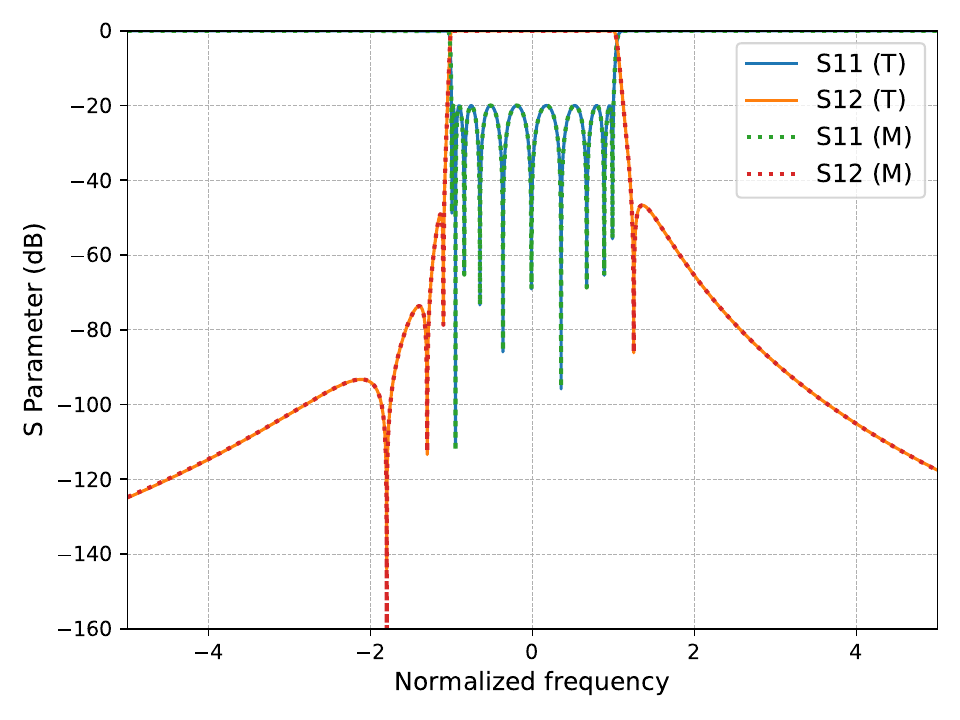}
           \caption{N=10,Folded}
       \end{subfigure}
       \begin{subfigure}[b]{0.24\linewidth}
           \centering
           \includegraphics[width=\linewidth]{k=5,folded.pdf}
           \caption{N=10,Extended Box}
       \end{subfigure}
       \caption{S-parameters of the target response generated by T and  M extracted by SMLP.}
       \label{res}
   \end{figure*}
   
\textbf{Analysis.}
From the data presented in Table \ref{Filter}, it is evident that the SMLP achieves a success rate of 100\% for target matrices for Folded topology, and neural networks can converge to the specified accuracy very quickly. 
For the more complex Extended Box topology, SMLP also achieves a success rate of 82\%, which is very high compared with the results in \citet{pfluger2015coupling}, demonstrating the efficacy of the proposed SMLP.
Additionally, analysis of Figure \ref{res} reveals that the response curves of matrices $T$ and $\textbf{M}$ are identical. This indicates that the coupling matrix $\textbf{M}$ obtained using the proposed SMLP method exhibits consistent responses with $T$, meeting practical requirements.

\section{Conclusion}
This paper is the first attempt on solving algebraic inverse eigenvalue problems using networks. A novel neural network named SMLP to deal with the structured inverse eigenvalue problems is proposed. One novelty is the derivation of a unified framework for several SIEPs, including the NIEP, SNIEP, StIEP, GStIEP and PEIEP. The other novelty lies in a hard constraint strategy, where the required orthogonal constraint is embedded into the network using an orthogonal decomposition of matrix.  
Several numerical examples from literature and engineering illustrate the efficiency of the proposed network SMLP. Besides, the effect of activation functions is also investigated numerically. The results indicate that the function $\tanh$ is more robust, especially for large scale problem, and the reason is still missing.
In addition, the idea of hard constraint strategy can be extended to many other problems where orthogonal constraint applied, such as the inverse singular value problems, computer vision, principal component analysis. 

\section*{Acknowledgements}
This work is supported in part by 
the National Key R\&D Program of China under the grant 2021YFA1003400, the National Natural Science Foundation of China under the grant 12071069,  
the Fundamental Research Funds for the Central Universities under the grant 2412022ZD032
and 
the Science and Technology Development Planning of Jilin Province under the grant YDZJ202201ZYTS573. 

\bibliographystyle{elsarticle-num-names}

\bibliography{cas-refs}

\end{document}